\documentclass[11pt, twoside, a4paper]{article}
\usepackage{amsfonts, amssymb, amsthm}
\usepackage[latin1]{inputenc}
\usepackage[intlimits,fleqn]{amsmath}
\usepackage[english]{babel}
\usepackage{graphicx}

\include{epsf}
\newcommand{\la}{\langle}
\newcommand{\ra}{\rangle}

\newcommand{\suli}[2]{\sum\limits_{#1}^{#2}}
\newcommand{\ili}[2]{\int\limits_{#1}^{#2}}

\newcommand{\wh}[1]{\widehat{#1}}

\def\a{\alpha}
\def\ap{\alpha'}
\def\ba{\bar{\alpha}}
\def\b{\beta}

\def\Nn{|\mathcal{C}_o|}
\def\Nh{\widehat{N}}

\def\Mb{L}

\def\Cn{\mathcal{C}_o}

\def\cp{\chi_{p}}

\def\Em{\E_\mu}

\def\Pn{\Pi_{\Lambda_n}}
\def\eh{\frac{1}{2}}

\def\be{\begin{eqnarray*}}
\def\bel{\begin{eqnarray}}
\def\ee{\end{eqnarray*}}
\def\eel{\end{eqnarray}}

\def\L{\mathcal{L}}
\def\Ln{\Lambda_n}

\def\R{\mathbb{R}}
\def\N{\mathbb{N}}
\def\E{\mathbb{E}}
\def\P{\mathbb{P}}

\def\Pgl{\mathbb{P}\_}
\def\C{\mathbb{C}}

\def\F{\mathcal{F}}

\def\Y{\mathcal{Y}}
\def\Z{\mathbb{Z}}
\def\X{\widehat{X}}

\def\e{\epsilon}

\def\d{\delta}

\def\id{\mathbb{I}}
\def\Sp{\textrm{Tr}}

\def\nn{\nonumber}

\def\eqiv{\Leftrightarrow}

\def\vec{\left(\begin{array}{cc} }
\def\cev{\end{array} \right)}
\def\w{\omega}

\theoremstyle{plain}
\newtheorem{theo}{Theorem}[section]
\newtheorem{lemma}[theo]{Lemma}

\newtheorem{corollary}[theo]{Corollary}

\theoremstyle{definition}
\newtheorem{defi}[theo]{Definition}

\begin{document}
\begin{titlepage}

{\center
{\large 

\quad  {\bf An interlacing technique for spectra of random walks and its
  application to finite percolation clusters}}

\vspace{0.5cm}

Florian Sobieczky \footnote{Institut f\"ur Mathematik C, TU-Graz, supported by
  FWF (Austrian Science Fund), project P18703.}\\ 
}

\vspace{0.5cm}

{\bf Abstract:} A comparison technique for random walks on finite graphs
is introduced, using the well-known interlacing method. It yields improved
return probability bounds. A key feature is the incorporation of parts of
the spectrum of the transition matrix other than just the principal eigenvalue. As an
application, an upper bound of the expected return probability of a
random walk with symmetric transition probabilities is found. In this case, the
state space is a random partial graph of a regular graph of bounded geometry
and transitive automorphism group.  The law of the random edge-set is assumed
to be stationary with respect to some transitive subgroup of the automorphism
group (`invariant percolation'). Given that this subgroup is unimodular, it is
shown that stationarity strengthens the upper bound of the expected return
probability, compared with standard bounds derived from the Cheeger inequality.\\

{\bf Keywords:} random walks, random walks on random partial graphs,
percolation, critical percolation, heat kernel decay, return probability,
comparison theorems\\ 

{\bf AMS classification:} 60K37, 60B99, 60D05, 60J10, 60J27

\tableofcontents

\end{titlepage}  
\newpage

\section{Interlacing for Random Walks on Graphs}

How does the return probability of a simple random walk (SRW) on a graph change
under removal and insertions of edges? The intention of the present paper is to
answer this question for another specific random walk on finite graphs. The
link to SRW will be provided in the form of inequalities comparing these two
Markov processes.  Using a circle of results from matrix theory called {\em
  interlacing}, we derive a method for comparing the return probabilities of
random walks of the same type, but realized on different graphs. Moreover, we
apply these results to invariant percolation with finite clusters. While the
first section of this paper is concerned with the presentation of the main
result in the context of standard interlacing theory, the second part describes
the implications concerning the annealed return probability on finite,
percolative subgraphs of transitive, unimodular graphs. The third section is
devoted to the proof of the main result and a geometric property of finite
trees.
\subsection{Introduction}

Partial graphs are subgraphs in which only edges are removed, while the set of
vertices remains the original one [\ref{cds}]. We will study a certain type of
random walk on partial graphs of graphs for which the return probability is
known. In the second section, we will consider finite partial graphs of {\em
  transitive} graphs. This will include allowing for connected components
consisting of only single vertices. The percolation generating these components
will be assumed only invariant. The finite clusters may occur as connected
components of either subcritical, supercritical, or critical percolations.\\

{\em Interlacing} refers to a set of techniques concerning the spectrum of
matrices under perturbations of known rank. The smaller the rank of the
perturbing operator, the more similar the set of eigenvalues of two matrices,
one being the perturbed version of the other. The heart of the interlacing
methods is the Courant-Fischer variational principle. In the present context,
it will be used to show a monotonicity property of the spectrum of a certain
regularised form of a given simple random walk, which we will call the {\em
  re- gularised random walk} (RRW): under removal of edges, the eigenvalues of
the transition matrix increase. Furthermore, insertion of edges results in a
shift of the $k$-th eigenvalue across less than $k$ intervals of the
unperturbed spectrum [\ref{hj}].\\

A similar property holds for the spectrum of combinatorial Laplacians of graphs
[\ref{haemers}, \ref{heuvel}]. The combinatorial Laplacian of subcritical Bernoulli-bond
percolation graphs has been considered in [\ref{kirsch}] for $\Z^d$, and
amenable graphs in [\ref{antves}]. For the adjacency matrix, as well as
`Laplacian spectra of graphs', interlacing is successfully applied (e.g.
[\ref{guo}, \ref{guo2}]).  In [\ref{chen_und_co}, \ref{li}] interlacing
results for normalised Laplacians have been proven. In the present paper,
interlacing will be used to find lower bounds for the eigenvalues of the
transition kernel of the RRW, which is, up to scaling, equivalent to the
Laplacian spectrum of graphs.  RRW has been considered (not under this name) in
many other contexts.  In [\ref{cheper2}] and [\ref{lyons}], RRW is called {\em delayed  random walk}. \\

For reversible random walks, it is a fundamental and well-known fact that the
spectral properties of the transition probability kernel are directly linked to
the quickness of the Markov chain to approach stationarity. In particular, the
{\em spectral gap} ([\ref{laurent}], sec. 2.1.2) as a function of the geometry
of the random walk's state space (given by a finite connected graph, here)
plays a characteristic role in the determination of the rate of convergence in
terms of the graphs order (=: size of its vertex set). If $\lambda$ denotes the
spectral gap, then it is well known that the return probability of a reversible
random walk $P_t(o,o)$ of some initial site $o$ of a connected graph with
$N<\infty$ vertices obeys ([\ref{chung}], chap.  1.5), with $c$ some constant
generally depending on $N$ and on the degree of $o$,
\bel 
|P_t(o,o) - \pi(o)|\;\;\le \;\;c\,e^{-t\lambda},\label{eq:1} 
\eel

where $\pi(\cdot)$ denotes the stationary distribution. A Characteristic of the
first eigenfunction of a discrete Laplacian with Neumann boundary conditions is
$\pi \equiv 1/N$. Under these circumstances $\lambda \ge c'/N^2$, for some
$c'>0$, by Cheeger's inequality (see [\ref{diastr}], Corr. 2.1.5).\\

Our goal in this paper are upper bounds for the annealed return probability of
the \mbox{regularised} and simple random walks on a random finite graph. To achieve
this, geometric information of the subgraph will be incorporated into spectral
estimates.  In particular, more information about the spectrum of the
transition kernel than just the spectral gap will be used. Instead of comparing
each of the eigenvalues with the principal eigenvalue, an estimate is derived
by using an interlacing technique.  It results in an estimate similar to that
in (\ref{eq:1}), given by Theorem \ref{theo:haupt}, however, for the return
probability $\bar{P_t}=(1/N\sum_kP_t(k,k))$. This is the case of a uniform initial
distribution (averaged return probability). For $t$ from a certain
`intermediate time-window', there is an additional polynomial prefactor: \bel
|\bar{P_t} - \pi(o)|\;\;\le\;\;a\frac{N^r}{t^{s}}e^{-t\lambda}.  \eel

Here $a$ is a known constant, and $r<2s$. This strict inequality is the reason
for an improvement over using $(1)$ to obtain bounds for the annealed return
probability. The constants $a,r$, and $s$ depend {\em only} on the largest occurring
degree, as is $c$ in (\ref{eq:1}) for the RRW. The consequences for the annealed
return probability of the RRW on finite clusters of invariant percolations with
exponential decay of the cluster-size are discussed in section
\ref{sec:invper}. It will be seen that the time-window spans far enough to
result in the aforementioned improvement.\\ 

As for notation, apart from a finite simple random walk $X_n, n\in \N$ in
discrete time with transition probability matrix $P$, we will consider the
corresponding regularised random walk $\widehat{X}_n, n\in\N$, with transition
probability matrix $A$. We have chosen this letter to remind the close
connection to the adjacency matrix of the regularised graph whose vertex set is
the state space. The adjacency matrix of the original graph is denoted by
$\mathcal{A}$. By $D$ we denote the diagonal matrix with entries given by the
vertices' degrees.  $\sigma(A)$ is the set of eigenvalues of $A$.  Since
we will need both terms, we distinguish between {\em combinatorial} Laplacian
$L=D-\mathcal{A}$ of a graph ({\em admittance}, [\ref{cds}]), and its {\em
  normalised} Laplacian $\mathcal{L} = \id - D^{-\eh}\mathcal{A}D^{-\eh}$
[\ref{chung}].  The Laplacian spectrum of a graph is $\sigma(L)$ (e.g.
[\ref{guo}]).  Finally, subtracting an edge from a graph will be written in
short form by $G-e$, which means that the vertex set remains unchanged, while
the edge set looses `$e$'.  We will reserve the letter $H$ for finite graphs,
while $G$ will denote infinite, transitive graphs. In particular, $H$ will be
the finite, random connected component containing the root of some $G$ with
bounded geometry of degree $\d$.

\subsection{Interlacing: standard results}\label{sec:interlacing}

In matrix theory, the term {\em interlacing} refers to the nature of the change
of the spectrum of a hermitian matrix under a hermitian perturbation of given
rank [\ref{hj}, \ref{lan}]. In particular, the location of the eigenvalues on
the real line of the perturbed matrix may be compared to the eigenvalues of the
unperturbed one. These techniques can be extended to the singular values of a
general matrix and are used in various contexts.\\

In graph theory, interlacing has been used to study the spectrum of the
adjacency matrix and the Laplacian of a finite graph, e.g. by [\ref{haemers},
\ref{heuvel}, \ref{godsil}, \ref{bolobas}, \ref{guo}], but also by many others
(see [\ref{haemers}] for a survey).  In [\ref{bro}], interlacing is used in
connection with the return time of random walks. In the context of random
walks, we like to show that by considering transformations of graphs which
connect different components, it is useful to perform comparisons of the return
probabilities of {\em reducible} random walks.

\begin{figure}[htb]
\centerline{
\epsfxsize=10cm
\epsffile{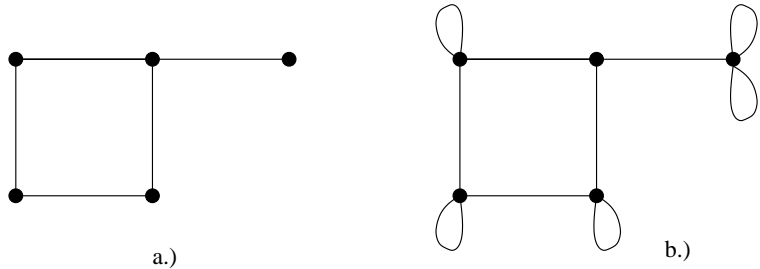} \hskip 0.0cm}
\caption{a.) Finite connected graph, with highest degree equal to three
  ($\d=3$); \newline \hfill b.) regularisation of the graph (RRW is the simple
  random walk on the regularised graph.)}\label{fig:reg_walk}
\end{figure}

\begin{defi}\label{defi:eins} Given a finite, simple (single-edged, undirected,
  loopless; see [\ref{doug}]) graph $H=\la V,E\ra$, let $\widehat{H}$ be its
  regularisation with loops, where for each loop attached to a vertex, its
  degree is increased by {\em one}. Let $(\widehat{X}_n, \mu_0)$ be the simple
  random walk on $\widehat{H}$ with initial distribution $\mu_0$ and call it
  {\bf regularised random walk} on $H$ (see Fig. \ref{fig:reg_walk}).
\end{defi}

Many of the properties of the regularised random walk are passed on to the
simple random walk, due to the following `equivalence' (see also
[\ref{cheper2}], paragraph 3 ).  

\begin{defi} The return probability of a random walk $X_n, (n\in\N)$ with
  finite state space and {\em uniform initial distribution} into its initial
  state $X_0$, is called {\bf average return probability}, denoted by
  $\Pgl[\wh{X}_n = \wh{X}_0]$ .
\end{defi}

\newpage

The following definitions are standard in the theory of finite random walks
[\ref{laurent}]  and spectral graph theory [\ref{chung}]:

\begin{defi}(Time-continuous random walk)
  The time-continuous regularised random walk $\wh{X}_t$ associated with the
  discrete-time regularised random walk $\wh{X}_n$ on a graph $H=\la V,
  E\ra$ with transition matrix $A$ is the random walk performing $n$ discrete
  independent steps across edges in $E$ within the time interval $[0,t]$ with
  probability $e^{-t}t^n/n!$. Therefore, in the case of  $H$ being finite, if
  $(e, A f) = \sum_{k,l \in V} e(k) A_{kl} f(l)$, the transition probability
  from $k$ to $l$ in finitely many steps is given by
\be
\P^k[\wh{X}_t = l]\;=\;\suli{n=0}{\infty} e^{-t}\frac{t^n}{n!}(e_k, A^n e_l) =
(e_k, e^{-t(\id-A)} e_l).
\ee

The quadratic form $(e_k, e^{-t(\id - A)} e_l)$ with $e_k,e_l \in \R^{|V|}$
the unit vectors along the direction corresponding to $k,l\in V$, respectively,
 is called the {\bf heat kernel} of $A$. It is the kernel of the
time-continuous semigroup $\exp(-t(\id - A)), (t\ge 0)$ associated with
$A$. The time-continuous simple random walk $X_t$, which is  associated with
$P$, is defined analogously, with $A$ replaced by $P$.
\end{defi}

We denote by $l^2(V, \pi)$ the space $\{\; f:V\to \C\;:\;(f,f)_\pi  < \infty
\;\}$, where $(f,g)_\pi:=\sum_k f(k)\overline{g(k)}\pi(k)$, and $\pi\in \mathcal{M}_{+,1}(V)$
is some probability distribution on  $V$. Naturally, for finite $V$, with
$|V|=N\in \N$,  this is isopmorphic to $\C^N$. 

\begin{defi}(Laplacians and heat kernel for the regularised graph) The graph
  Laplacian of the regularised graph $\wh{H}$ is denoted by
  $\wh{L}=\wh{D}\;-\;\wh{\mathcal{A}}$, where $\wh{D}=\d\id$, and the
  adjacency matrix $\wh{\mathcal{A}}$ accounts for the attached loops:
  $\wh{\mathcal{A}}\;=\;\mathcal{A}\;+\;\d\id\;-\;D$. Therefore,
\be
L \;=\; D \;-\; A \;=\; \wh{D}\;-\;\wh{\mathcal{A}}\;=\;\wh{L},
\ee

i.e. the graph Laplacian $L$ is invariant under the decoration of $H$ with
loops. $D=\textrm{diag}(\textrm{deg}(v_1), \textrm{deg}(v_2),\cdot\cdot\cdot,
\textrm{deg}(v_N))$ and deg($v_j$) the degree of vertex $v_j$ of $H$. The
normalised Laplacian for the regularised $\wh{H}$ is $\wh{\L}= \wh{D}^{-\eh}
\wh{L} \wh{D}^{-\eh}=(1/\d)L$.  The quadratic forms $(e, \L f)_\pi$ and
$(\wh{e}, \wh{\L} \wh{f})_{\wh{\pi}})$ with $f,e\in l^2(V, \pi)$ and $\wh{f},
\wh{e}\in l^2(V, \wh{\pi})$ where $\pi$ and $\wh{\pi}$ are the invariant
distributions of $X_t$ and $\wh{X}_t$ are called {\bf Dirichlet-Forms} of $\L$
and $\wh{\L}$, respectively.
\end{defi}

We quote the following  comparison theorem for the eigenvalues of the random
walks from [\ref{laurent}], section 1.2.5 without proof in a simplified form. It is used in the
proof of the next lemma.

\begin{lemma}\label{lemma:comp_dir} Let $\L$ and $\L'$ be two
  normalised Laplacians on $l^2(V,\pi)$ and $ l^2(V, \pi')$, respectively.
   Let there be a linear map $l^2(V,\pi)\ni f\mapsto f'\in l^2(V,\pi')$ and
   $a>0$ such that for every $f\in l^2(V,\pi) $ it holds   that
$$ (f, \L f )_\pi \;=\;  (f',  \L' f' )_{\pi'}\;\;\;\;\;\;\;\textrm{ and }\;\;\;\;\;\;\;a\cdot(f,f)_\pi \le (f',f')_{\pi'}. $$

Then it holds for all pairs of $j$-th eigenvlaues $\lambda_j$ and $\lambda'_j$
of $\L$ and $\L'$, repectively,
$$ a\cdot\lambda_j \le \lambda'_j, $$
 if both sets of eigenvalues ( $\sigma(\L)$ and $\sigma(\L')$) are enumerated
 in the same way (e.g. $\lambda_j\le \lambda_{j+1}$ and $\lambda'_j\le
 \lambda'_{j+1}$ for all $j\in\{1, ...., |V|-1\}$).
\end{lemma}

\newpage

{\bf Remark:} The average return probability of the time-continuous $X_t$ and
$\wh{X}_t$ relates to the trace of the corresponding heat kernel (see
[\ref{laurent}, \ref{mathieu}]). While $\exp(-t(\id - P))$ is the heat kernel of $X_t$, for
the regularised random walk it is given by $\exp(-t(\id - A))$,
where $A=\d^{-1}\wh{\mathcal{A}}$, with $\wh{\mathcal{A}}$ the adjacency matrix
of the regularised graph $\wh{H}$. Then, with $|V|=N$:
\be 
  \Pgl[\wh{X}_t =  \wh{X}_0] \;=\;\frac{1}{N} \suli{k\in V}{} \P^k[\wh{X}_t=k] \;\;=
  \;\; \frac{1}{N}\Sp \exp(-t(\id-A)).  
\ee

\begin{lemma}\label{lemma:spurmono}
  Assume that the finite graph $H=\la V,E\ra$ is connected. Then the average continuous time
  return probabilities of $X_t$ and $\wh{X}_t$ on $H$ fulfil
\be
\Pgl[\;X_t\;=\;X_0\;]\;\;\le\;\;\
\Pgl[\;\wh{X}_t\;\;=\;\;\wh{X}_0\;]\;\;\;
\le\;\;\Pgl[\;{X}_{t/\d}\;=\;{X}_0\;].
\ee
\end{lemma}

{\em Proof:} Introduce $\bar{P}=(1-1/{\d})\id+(1/{\d})P$, the transition-matrix
of the {\em lazy} random walk, which is the simple random walk on the graph
$\bar{H} = \la V, \bar{E}, \bar{m}\ra$, with $N=|V|$ vertices, and edges
$\bar{E}=E\cup \bar{E}'$, where $\bar{E}'$ are the multiple loops
$\{k\}\in\bar{E}$ with multiplicity 
\be 
m(\{k\}) = ({\d}-1) \cdot\textrm{deg}_H(k).
\ee

Then the lazy normalized Laplacian is
$\bar{\L}=1-\bar{D}^{1/2}\bar{P}\bar{D}^{-1/2}=(1/\d)\L$, a scaled version of
the normalized Laplacian of $H$, where $\bar{D}=D+({\d}-1)D=\d D=\textrm{diag}(\textrm{deg}(v_1), \textrm{deg}(v_2),\cdot\cdot\cdot\cdot,
\textrm{deg}(v_N))$ is the diagonal degree operator of $\bar{H}$ taking account
of the loops multiplicity.  Note that
the combinatorial Laplacian $L=D-\mathcal{A}$ (or admittance operator) which
doesn't change under the regularisation (s.a.), also remains invariant under
`laziness', i.e.
\bel 
D\;-\;\mathcal{A} \;=\;{\d}\id\;-\;\hat{\mathcal{A}} \;=\; \bar{D} \;-\; \bar{\mathcal{A}}.\label{eq:Ls} 
\eel

Here, $\mathcal{A}$, $\wh{\mathcal{A}}$, and $\bar{\mathcal{A}}$ are the adjacency
matrices of $H,\wh{H}$ and $\bar{H}$ taking into account the changed degree due
to the loops. ${\d}\id$ and $\bar{D}=D+({\d}-1)D=\d D$ are the corresponding `degree
operators' for $\wh{H}$ and $\bar{H}$. \\

In this situation, the properties of the normalised Laplacian differ from
those of the combinatorial Laplacian. Instead of invariance under
`loop-decoration', the spectrum $\{\bar{\lambda}_j(\bar{\mathcal{L}})\}$ of the
lazy $\bar{\mathcal{L}}$ is a scaled version of $\sigma(\mathcal{L})=
\{\lambda_j(\mathcal{L})\}$: for $j\in\{1, \cdot\cdot\cdot, N\}$,
\bel 
\lambda_j(\bar{\mathcal{L}})\;=\;\frac{1}{\d}\lambda_j(\mathcal{L}).
\label{eq:lambda_eq} 
\eel

On the other hand, the diagonal elements of $D$ fulfil $(v,D v)\le(v,(\d \id)v)\le
(v,({\d}D) v)$. Therefore, while the Dirichlet-forms (for $L, \wh{L}$, and
$\bar{L}={\d}D-\bar{\mathcal{A}}$) are equal in value, the norms
$\|\cdot\|_{\pi}$ of $v$ in $l^2(V, \pi)$, where $\pi$ is the invariant
distribution of $P,A$, or $\bar{P}$, differ: $\|v\|^2_\pi=(v,v)_\pi=(v,Dv)\le
(v,\wh{D}v)=(v,v)_{\wh{\pi}}\le (v, \bar{D}v)=(v,v)_{\bar{\pi}}$ for $v\in l^2(V,\pi)$.  By Lemma 
\ref{lemma:comp_dir}, it follows that  
\bel 
\lambda_j(\mathcal{L}) \;\;\;\ge\;\;\; \lambda_j(\wh{\mathcal{L}})
\;\;\;\ge \;\;\;\lambda_j(\bar{\mathcal{L}}).\label{eq:compare_lambda}
\eel

Therefore, $\suli{k\in V}{} \P^k[X_t = k] = \Sp[\exp(-t\mathcal{L})] \le
\Sp[\exp(-t\mathcal{\wh{L}})] \le \Sp[\exp(-t\mathcal{\bar{L}})]. 
$\hfill\qed\\

For completeness, we add the following standard result about the relevant
comparison between discrete time random walks:

\begin{lemma} \label{lemma:comp2}
The average return probabilities for discrete ($\P^d_\_[\cdot]$) and
  continuous ($\P^t_\_[\cdot]$) time on $H$ of $(X_\cdot, \mu_0)$ and
 (with $\mu_0$ the uniform distribution on $V$) fulfil\\[-0.5cm]
\be 
  \Pgl^d[X_{4n} = X_0] \;\;\; \le\;\;\; 2\;\Pgl^t[X_{2n} = X_0 ]\;\;\; \le\;\;\;
  2\;\Pgl^d[X_{n} = X_0 ]+2e^{-t}.   
\ee
\end{lemma}

{\em Proof:} This follows from standard methods of comparison between discrete,
lazy and continuous-time versions of a given Markov chain 
(see [\ref{woess}], Lemma 14.2.c [\ref{laurent}], and Corollary 1.3.4): Let
$\bar{P}=\eh\id+\eh P$.  With $\Sp[\exp(-t(1-P))]=\Sp[\exp(-t\L)]$ it holds due to $(1-x)^n\le\exp(-nx)$ for $0\le x\le1$, and
$\log(1-y)\ge -2y $ for $0\le y \le \eh$  
\be 
\Sp[P^{4n}]\;\le\; 2\;\Sp[\bar{P}^{4n}]
\;\le\;
2\Sp[e^{-4n(1-\bar{P})}]\;\;=\;2\;\Sp[e^{-2n\mathcal{L}}]\;\le\;2(\Sp[P^n]+e^{-t}). \;\;\;\;\;\;\; \hfill \qed
\ee

From now on let $\sigma(A)=\{ \beta_j \}_{j=1}^N$ be the
eigenvalues of $A$ enumerated in a non-increasing way, i.e. $1=\beta_1 > \beta_2
\ge \beta_3 \ge \beta_4 ... \ge \beta_N$. The strict inequality arises from
$H$ being connected, with $\beta_1$ being the Perron-Frobenius eigenvalue of
$A$. The restriction of the interlacing method to the transition kernel of the
{\em continuous} time RRW, for which the eigenvalues $\exp(-t(1-\beta_j))$ are
all positive, and therefore easier to handle in comparison techniques, is now
feasible.  The bounds obtained for this setting can then be transferred to
discrete time, and the simple random walk, by the use of the two Lemmata
\ref{lemma:spurmono} and \ref{lemma:comp2}.\\

\begin{figure}[htb]
\centerline{
\epsfxsize=10cm
\epsffile{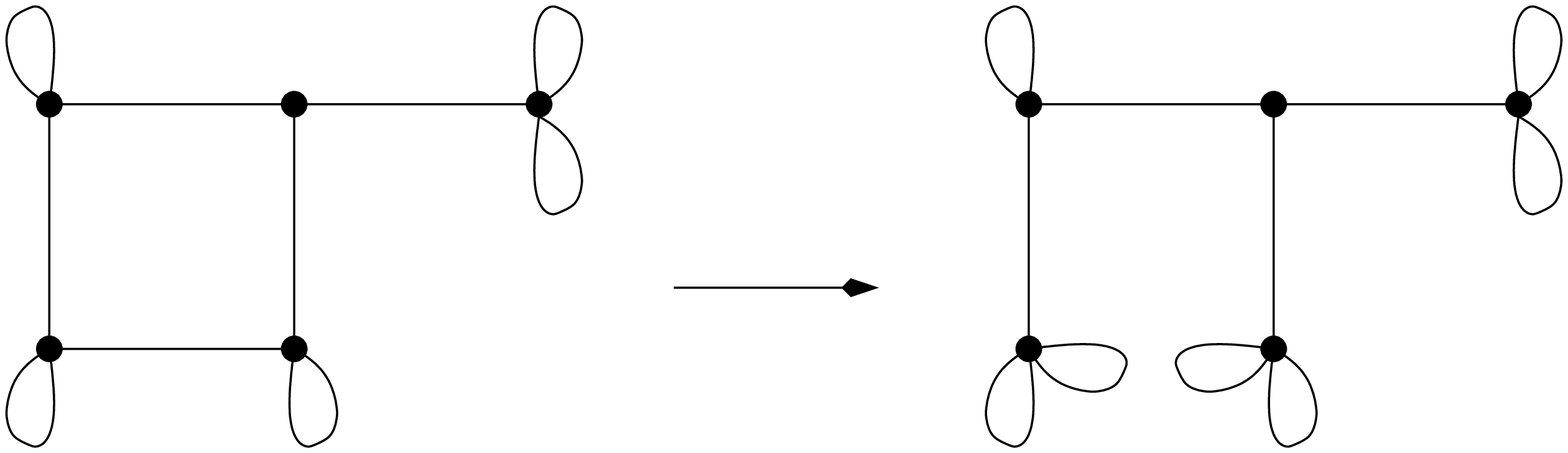} \hskip 0.0cm}
\caption{Removing an edge of a (loopless) graph acting as the state space of
  RRW corresponds to the SRW on a decorated graph with the edge {\em replaced}
  by a pair of loops.}\label{fig:trans_walk}
\end{figure}

Now, we find bounds for eigenvalues of the RRW on finite graphs with edges {\em
  removed}: Let $A'$ denote the transition probability matrix of the
regularised random walk of the graph $H'=H-e$, obtained by removing a single
edge and keeping all vertices (see Fig. \ref{fig:trans_walk}).

We then have $A' = A + S_e$, with $S_e$ of the following type:
\bel
S_e\;=\;\frac{1}{\d}\left(\begin{array}{ccccc}
    0& & ...& &0\\
    ...&+1& ...&-1&...\\
    ...& ...&...&...&...\\
    ...&-1& ...&+1 &...\\
    0& & ...& &0 \end{array} \right),\label{eq:matrix}
\eel

where the non-zero elements appear in the rows and columns indexed by the
vertices of the edge $e$ removed. This is $(1/\d)$ times a projection onto the
one-dimensional subspace spanned by $(0,...,0, 1,0, ....,0,, -1,0, ...,0)$, where the
non-zero entries occur at the positions indexed by the vertices incident with
$e$. It follows that $S_e$ is positive semi-definite, so $\beta_j\le \beta_j'$.  Moreover, its
rank equals one. Taking away all edges from some subset $R\subset E$, with
$r=|R|$, leads to
\bel
A' = A + \suli{e \in R}{} S_e.\label{eq:perturbation}
\eel

This is a perturbation with rank bounded by $r$. (Note, if $R$ are the edges
removed from the initial graph, the number of linearly
independent conditions induced by the boundary conditions may in general be smaller than $r$.)
 By the Courant-Fisher variational principle (see [\ref{hj}, \ref{lan}]), we have, with $\Y$
a subspace of $\C^N$ and $S=\suli{e \in R}{} S_e$, if $j>r$,

\bel
\beta'_{j} &=& \max\limits_{\Y \subset \C^N \atop \dim(\Y) = j} 
\;\;\min\limits_{v\in\Y} \frac{(v,A' v)}{(v,v)}\nn\\
&\le& \max\limits_{\Y \subset \C^N \atop \dim(\Y) = j} \;\;\;\; \min\limits_{ v\in\Y:
  \; \; { v_k=v_l } \atop {f.a. \;\{k,l\}\in R  }}
\;\;\frac{(v,\;(A\;+\;S \;) \;v)}{(v,v)} \nn\\ 
&=& \max\limits_{\Y \subset \C^N \atop \dim(\Y) = j} \;\;\;\; \min\limits_{ v\in\Y: 
  \; \; { v_k=v_l } \atop {f.a. \;\{k,l\}\in R  }}
\;\;\frac{(v,\;A \;v)}{(v,v)}\nn\\
&\le& \max\limits_{\Y \subset \C^N \atop \dim(\Y) =
  j-r}\min\limits_{v\in\Y} \frac{(v,A v)}{(v,v)} \;\;=\;\;
\beta_{j-r}. \label{eq:cf}
\eel
If $j\le r$, since $A$ and $A'$ are stochastic matrices, $\beta'_j \le
\beta_1=1$. \\[-0.2cm]

\begin{theo}(Interlacing for RRW on finite partial graphs)\label{theo:partial}
  The eigenvalues $\{\beta'_j\}$ of a regularised random walk on a partial
  graph $ H'=H-R$ relate to those ($\{\beta_j\}$) of the regularised random
  walk on the initial graph $H$, with $r=|R|$, by 
\bel \beta_j
  \;\;\le\;\;\beta'_j\;\;\le\;\;\beta_{(j-r)\vee 1}. \label{eq:il_ineqs}
\eel
\end{theo}
{\em Proof:} Removing the edges in $R$ is a perturbation of type
(\ref{eq:perturbation}) and $r$ is an upper bound for its rank. The result
follows from the positivity of $S_e$, (\ref{eq:cf}), and by considering that
stochasticity of the matrix is preserved. \hfill \qed\\

The second inequality in (\ref{eq:il_ineqs}) gives an upper bound for
the change of eigenvalues of the RRW under {\em insertion} of edges.\\

Similar results have been proven for normalised Laplacians of non-regular
graphs [\ref{chen_und_co}, \ref{li}]. In this case the monotonicity is lost due
to the loss of semi-definiteness of the corresponding pertubation matrix
(instead of $S_e$), which, however, is still of rank one.\\

{\bf Remark:} The condition (constraint) given by $v_k=v_l$  is the same as
the vanishing discrete gradient `d' of $v$ as a function of the vertices (if some
direction of the edge $e=\{k, l\}$ has been agreed upon): 
\be v_k - v_l\;\;=\;\;\textrm{d}v (e) \;\;=\;\; 0.  \ee 

In this form, the constraint may be considered a Neumann boundary condition for
the corresponding normalised Laplacian $\L=\id \;- \;A$ (comp. with
[\ref{barry}]). This is due to the monotonicity of the shifts in the spectrum,
under an increase of `boundary', which is negative, similar to Laplacians of
functions in continuous function spaces. Note that increasing the boundary is
realised in the discrete setting by the removal of edges.  However, in the
context of {\em induced subgraphs} ([\ref{chung}], chapter 8.3), the Neumann
eigenvalues of a graph relate to yet another random walk, which shows a
different behaviour at `corners' - owing to a special definition of the
discrete gradient at 'irregular' parts of the boundary. It is easily shown
that the Neumann random walk on induced subgraphs also shows the monotonicity
of the spectral shift under removal of edges. In the present context, however,
the perturbed graphs are {\em partial graphs}, i.e. subgraphs of which only
edges (not vertices) are removed.\\

Before stating the main result, we prove a simple well-known interlacing
theorem for the Laplacian spectrum of partial graphs
[\ref{heuvel}, \ref{chen_und_co}]:

\begin{lemma}\label{lemma:li}

  Let $H=\la V, E\ra$ be a finite, possibly unconnected, undirected
  graph with $N=|V|$ vertices. Let $D=\textrm{diag}(..., $deg$(v_j), ...)$ be the
  $N\times N$ diagonal matrix with the degrees of the vertices as its entries,
  and $\mathcal{A}$ the adjacency matrix of $H$. Let $L=D-\mathcal{A}$ denote
  the combinatorial Laplacian of $H$ and $\sigma(L)=\{
  \;\lambda_j\;\}_{j=1}^{N}$ its spectrum. Let the eigenvalues be enumerated in
  an increasing way:
  $0=\lambda_1 \le \lambda_2 \le  \cdot\cdot\cdot \le \lambda_N$.\\

  Then, removing a single edge results in a graph with eigenvalues shifted
  downwards. In other words, in a transformation $H\to H'=H-e$, for some edge
  $e\in E$ and corresponding transformation $L \to L'$ of the combinatorial
  Laplacian, the $j$-th eigenvalue in $\sigma(L')$ will obey (with $j\in\{2,
  \cdot\cdot\cdot, N\}$):
\bel
\lambda_{j-1} \;\le\; \lambda_j' \;\le\; \lambda_j.\label{eq:li}
\eel
\end{lemma}

{\em Proof:} We note that decorating a graph with loops leaves the
combinatorial Laplacian invariant. Therefore, if $\d$ is the maximum degree of
$H$, regularising $H$ by adding loops to each vertex until every vertex has
degree equal to $\d$, where the degree of a loop counts one, yields the same
Laplacian. Call the graph obtained by the loop-decoration $\widehat{H}$ and its
Laplacian $\widehat{L}$.  So, $L=\widehat{L}$, and, in particular, the $j$-th
eigenvalue of $L$ equals the $j$-th eigenvalue of $\widehat{L}$.\\

Removing an edge $e$ of $H$ yields the transformed graph $H'=H-e$. Let the
corresponding decorated (regularised) graph be called $\widehat{H'}$ and its
Laplacian $\widehat{L'}$ (again, the hat denotes the decoration with loops): Then the
matrix $\widehat{L'} - \widehat{L}$ is negative semidefinite and of rank one,
since it is of the form (\ref{eq:matrix}), multiplied by $-\d$. Since
$\widehat{L}$ and $\widehat{L'}$ are real and symmetric, it follows from
$\wh{L}=\d\id \;-\;\d A$ with $A$ the transition matrix of the RRW that the
$j$-th eigenvalues $\widehat{\lambda}_j$ and $\widehat{\lambda'}_j$ of
$\widehat{L'}$ fulfil $\wh{\lambda}_j=\d(1-\beta_j)$ and
$\wh{\lambda}_j'=\d(1-\beta_j')$. Therefore, it follows from Theorem
\ref{theo:partial} for $j\ge 2$,
\be \lambda_{j-1} \;\le\; \widehat{\lambda_j'} \;\le\; \lambda_j.  \ee

Since $\widehat{L'} = L'$, it follows $\widehat{\lambda_j'}=\lambda'_j$, the
$j$-th eigenvalues of $L'$.  \hfill \qed

\subsection{Main result}

The second inequality in (\ref{eq:il_ineqs}) gives upper bounds for eigenvalues
of a graph into which edges have been inserted - if $\beta'_j$ is viewed as the
`unperturbed eigenvalue'.  This is the main ingredient of the following
theorem. Its proof, given in the last section, uses Theorem \ref{theo:higher}.
Before stating the result we define some numerical constants, some of which
find repeated use in several of the theorems.  Again, given the finite, simple
graph $H$ (see Definition \ref{defi:eins}), $N=|V(H)|$ is the size of its
vertex set, $\d$ is the maximum degree.  Furthermore, we define the following
constants; they will be given again at the places where they are
used: 
\bel
\nu:= \frac{1}{(\d-1)\ln 16 }, \;\;\;\;\; a
=1\;+\;\frac{\sqrt{\pi\d}}{2^{1-\nu}}, \;\;\;\;\;\check{t}=4^{1+\nu}N^{1-\nu},
\;\;\;\;  \hat{t}=4^{-\frac{(1+\nu)^2}{\nu}}N^4.\label{eq:consts}\\\nn
\eel

\begin{theo}\label{theo:haupt}
  The return probability of the time-continuous RRW $\X_t,$ for  $t\in[\check{t}, \hat{t}]$
  on a finite, connected, simple graph $H$ of order $N$ and maximum degree
  $\d\ge 3$ obeys the bound

\bel \Pgl[\X_t\;=\;\X_0]
  \;\;\le\;\;\frac{1}{N}
  \;\;+\;\;a\left(\frac{N^{1-\nu}}{\sqrt{t}}\right)^\frac{1}{1+\nu} e^{-4 t/(\d
    N^2)}.\label{eq:haupt}\\ 
  \eel
\end{theo}

{\bf Remark:} For times $N^{2(1-\nu)}<t< N^2$, the exponential factor is still
close to one, while the prefactor may already be small. The
result is to be compared with the trivial bound 

\bel
\Pgl[\X_t\;=\;\X_0]
\;\;=\;\;\frac{1}{N}\Sp\left[ \exp -t(\id-A)\right]
\;\;\le\;\;\frac{1}{N}+\frac{N-1}{N}e^{-\frac{4t}{\d N^2}}, \label{eq:trivial}
\eel

resulting from $\exp(-t(1-\beta_2))\le \exp\{-4t/(\d N^2)\}$. For times as
large as $\hat{t}=O(N^4)$, the improvement of the prefactor is negligible in
comparison with the effect of the exponential factor. This justifies picking
$t\le\hat{t}$ from an `intermediate time-frame'.\\

We now derive a theorem on the eigenvalues of a regularised random walk, which
is proved by an interlacing technique involving {\em unconnected graphs}. It
will be used in the proof of Theorem \ref{theo:haupt}. It is easy to show that
the second largest eigenvalue of the transition kernel of RRW on any finite,
connected tree is bounded from above by the second largest eigenvalue of the
{\em path} (=: $P_N$) with the same number of vertices $N$.  With the spectrum
of the RRW on $P_N$ given by $\beta^o_j= 1\;-\;\frac{2}{\d}
(1\;-\;\cos\frac{(j-1) \pi}{N}), \;j\in\{1, ..., N\}$ (see [\ref{feller}], or
[\ref{laurent}], Ex.  2.1.1), this gives the upper bound for its second largest
eigenvalue $\beta_2^T$ (`T' for tree):
\bel 
\beta^T_2 \;\;\;\le\;\;\; \beta_2^o \;\;\;\le\;\;\; 1\;-\;\frac{4}{\d N^2}.
\label{eq:gen_second} 
\eel 

Using interlacing, we now find an upper bound for the eigenvalues of the
transition matrix of the RRW on a finite graph. For the proof, a lemma about
the geometry of finite trees is used, which is given in the last section.

\begin{theo}\label{theo:higher} Let  $\sigma(A)=\{\beta_j\}$ be the spectrum of the
  transition kernel of the regularised random walk on a finite, connected graph
  $H$ with $N$ vertices, enumerated in a decreasing way. Let $H$ have largest
  degree $\d$. Then, with $B\in\{0, ..., N-2\}$, and $\nu$ given in (\ref{eq:consts})

\bel 
\beta_{B+2}\;\;\le \;\; 1-\frac{4 }{\d N^2} \left(\frac{B+1}{2}\right)^{2
  \nu}.  \label{eq:higher_eq} 
\eel
\end{theo}

Before proving this theorem, we make some general observations and give
definitions regarding the interlacing inequalities (\ref{eq:il_ineqs}) applied
to arbitrary spanning trees of $H$:\\

Removing an edge shifts the spectrum of the RRW (of the transition probability
matrix) {\em upwards}, by at most one interval between successive eigenvalues.
Therefore, if we look at the $j$-th eigenvalue of the RRW on $H$ ($\beta_j$),
it is possible to bound $\beta_j$ from above by the corresponding eigenvalue
$\beta_j^T$ of RRW on any spanning tree. Likewise, it is possible to bound
$\beta_j^T$ from above by the $j$-th eigenvalue of the RRW on any forest
resulting from removal of edges: In the case where removal of edges results in
several connected components, the multiplicity of the eigenvalue equal to one
is equal to the number of connected components.\\

Comparing the second largest eigenvalue $\beta^T_2$ with an eigenvalue
$\beta_2^{T,1}$ of a graph resulting from removal of one edge (yielding a
disconnected graph with two trees as connected components) only yields
$\beta^T_2 \le 1$. The reason is that the multiplicity of `$1$' in the spectrum
of the disconnected graph is two and $\beta_2^{T,1}=1$. Define $\beta_j^{T,k}$
to be the eigenvalue of RRW on the spanning forest resulting from removal of
any $k$ distinct edges from an arbitrary spanning tree of $H$. Let $j\mapsto
\b^{T,k}_j$ be enumerated in decreasing order. Then, if the $j-$th
eigenvalue $\beta^T_j$ of RRW on the spanning tree of $H$ for $j\in \{2, \cdot\cdot\cdot, k+1\}$ is
compared with the $j$-th eigenvalue $\beta_j^{T, k}$, interlacing
(\ref{eq:il_ineqs}) yields for $j\in \{2, ..., k+1 \}$ only $\beta^T_j \le
1=\beta_1^{T,k}=\beta_2^{T,k}= \cdot\cdot\cdot \beta_{k+1}^{T,k}.$ \\

\begin{lemma}\label{lemma:precisely}  Removing $k$ edges from any spanning  tree of $H$, the
\mbox{$k+2$-nd} eigenvalue $\beta_{k+2}^{T,k}$ of the RRW on the resulting spanning
forest is known to be strictly smaller than one. In particular, it is the first
eigenvalue with this property in this enumeration:
\bel
\min\left\{\;j\in\{1, ..., N\}\;:\; \b^{T,k}_{j} < 1\;\right\} = k+2.\label{eq:min_j}
\eel
\end{lemma}

{\em Proof:} (Lemma \ref{lemma:precisely}) The forest resulting from removing $k$ edges from a
spanning tree of $H$ consists of precisely $k+1$ connected components.
Therefore, the eigenvalue with numerical value equal to one has multiplicity $k+1$.
 \hfill \qed\\

{\em Proof:} (Theorem \ref{theo:higher}) Let $\hat{\beta}_2^o$ be the second
largest eigenvalue of RRW on the path $P_{\hat{N}}$, decorated with loops to
make it a regular graph of degree $\d$, where $\hat{N}$ is the size of the
largest tree in the forest resulting from removing any $k$ edges of any
spanning tree of $H$. 
\bel
\beta_{k+2}\;\;\;\le\;\;\; \beta_{k+2}^T\;\;\; \le\;\;\; \beta_{k+2}^{T,k}\;\;\;
\le\;\;\; \hat{\beta}_2^o\;\;\;\le\;\;\;1\;-\;\frac{4}{\d \hat{N}^2}. 
\label{eq:rem_edges}
\eel

The first and second inequality is the first of (\ref{eq:il_ineqs}). By Lemma
\ref{lemma:precisely}, 
the third and fourth inequality of (\ref{eq:rem_edges}) are
(\ref{eq:gen_second}) with the order ($=N$) of the tree replaced by $\hat{N}$.\\

An optimal bound of $\hat{N}$ will be made by the help of Lemma
\ref{lemma:divide_trees} (compare with [\ref{guo}], where eigenvalues of the
adjacency matrix of a graph {\em without} loops are estimated in a similar way
- also involving interlacing).  We take an arbitrary spanning tree of $H$, and
note that we can choose a sequence of edges, which are taken away, such that
at each edge, the size reduces by at least a factor of $q=m/(1+m)$, where \be q
\;\;=\;\; 1\;\;-\;\;\frac{1}{4(\d-1)}.  \ee

Now, the order $\hat{N}$ of the largest remaining subtree depends on $B$, the
number of removed edges: If some $b\in \N$ is such that $2^b \le B+1 <
2^{b+1}$ (i.e. $b=[\log_2 (B+1)]$, the largest integer below $\log_2(B+1)$),
then by removing $B$ edges, the cardinality $N$ of the vertex set of $H$ may
be reduced by a factor of $q$ at least $b$ times. Therefore, since $q<1$,
\bel 
q^b \le q^{\log_2(B+1)-1} = \left(\frac{B+1}{2}\right)^{\log_2 q} =
\left( \frac{B+1}{2}\right)^{\log_2(1-1/(4(\d-1)))}, 
\eel

and $\hat{N} \le N q^{\log_2(B+1)-1} = N \left(\frac{B+1}{2}\right)^{\log_2 q} \le N/\left(\frac{B+1}{2}\right)^{\nu}$.\\

By (\ref{eq:rem_edges}), this implies 
$\beta_{B+2}\;\;\le\;\; 1\;\;-\;\;\frac{4 }{\d N^2}\left(\frac{B+1}{2}\right)^{2\nu}$.\hfill \qed \\

The proof yields another form of this theorem, namely, the direct comparison of the
spectrum of the RRW on $H$ by the spectrum of RRW on the path $P_N$:

\begin{corollary} \label{cor:higher} Under the conditions of the last
  theorem, let $\{\beta_j^o\}_{j=1}^{N}$ be the spectrum (ordered in
  decreasing order) of the transition kernel RRW on the path $P_{N}$, with a
  number of vertices $N$.  Then, with $x\mapsto [ x] := \max\{m\in \N\;:\;
  m\le x\}$, the eigenvalues $\{\beta_j\}$ of the RRW on $H$ relate to
  $\sigma(P_N)$ by \\[-0.2cm]
 \bel 
\beta_{B+1} \;\;\;\le\;\;\; \beta_{[ (B/2)^\nu ] + 1}^o,\;\;\;\;\;\;\;\;
  \;\;\;\;(B\in\{0, ..., N-2\}).\label{eq:higher}  \eel\\[-0.5cm]
\end{corollary}

{\em Proof:} Follows directly from the first three inequalities of
(\ref{eq:rem_edges}), and from

\be
\hat{\beta}^o_2 \;=\;  1\;-\;\frac{2}{\d}(1-\cos\pi\frac{1}{[ N (\frac{2}{(B+1)})^{\nu} ]})
\;\le\; 1\;-\;\frac{2}{\d}(1-\cos\pi\frac{[ 2^{-\nu}(B+1)^\nu]}{N})
\;=\;\beta^o_{[ 2^{-\nu}(B+1)^\nu] + 1}
\ee 

where $\hat{\beta}^o_2$ is the second largest eigenvalue of the path
$P_{\hat{N}}$, with $\hat{N}=[ 2^\nu N/(B+1)^\nu]$. \hfill\qed\\

{\bf Remark:} This is a comparison theorem for other elements of the spectrum
of $A$ than $\beta_2$ and, as such, compares to results concerning 'higher
eigenvalues' of the spectrum of the Laplacian $1-A$. (See e.g. [\ref{laurent}],
Theorem 3.3.17, where assumptions about the isoperimetric properties of the
graph are made.) The relevant geometric property assumed here is only the
uniform bound on the degree of the graph, given by $\d$.

\newpage

\section{Invariant Percolation}\label{sec:invper}

\subsection{The mass transport principle}

We apply our results to finite random graphs, by considering the {\em expected}
return probability (= annealed return probability) of a continuous-time random
walk on the finite random partial graphs of a transitive graph of finite
degree: namely the finite connected components of an invariant percolation.
All parameters and constants of the estimates shall be accessible, i.e. it
should be possible to express them as functionals (expected values) of the
distribution of the random process generating the subgraphs. Our results can be
used in this case, if the random graphs considered here will be restricted to
the subgraphs induced by the almost surely {\em finite} connected components.
subcritical percolation on the Cayley graphs of finitely generated groups are
included as a special case. See [\ref{kirsch}, \ref{antves}] for estimates of
the integrated density of states (`cumulative spectral measure') near the edges
of the generator's spectrum in this case.  However, also the critical
[\ref{BLPS}], or supercritical percolation measure, conditioned on the
finiteness of the cluster containing the root is a possible setting for the
results proven, here. See Theorems \ref{theo:2dcritical} and
\ref{theo:treecritical} for this subject. \\

One of the most important results in percolation theory in the last few years
was the fact that critical percolation on non-amenable, unimodular graphs has
almost surely finite clusters [\ref{BLPS}]. The {\em mass transport principle}
[\ref{schon}] states that for invariant percolation on an infinite, unimodular,
transitive graph $G=\la V, E\ra$, a function of two vertex-valued arguments,
invariant under diagonal action, allows interchange of the arguments under
summation of one of them, i.e. for all $f:V\times V\to\R$,
\be
\suli{v\in V}{} f(v,w)\;=\;\suli{v\in V}{} f(w,v).
\ee

See [\ref{lyons}] for a detailed discussion of the consequences in the theory
of invariant percolation. 

\subsection{Application to the expected return  probability} 

The results concerning the return probability of the RRW were formulated for an
initial distribution given by the uniform distribution on the finite cluster.
The relation to the return probability to a fixed, `deterministic' starting
point in the case of invariant percolation is settled for random partial graphs
of unimodular graphs by the following lemma. Consider $H(\w), \w \in \Omega$ be
a random partial graph of a unimodular transitive graph $G=\la V, E\ra$ realised
with probability given by an invariant law $\mu$ on the probability space
$\Omega=2^E$, where $\F$ is the product $\sigma$-algebra on $\Omega$. Let
$\Cn$, depending on $\w\in\Omega$, be the connected component of $H(\w)$
containing the root $o\in V$.

\begin{lemma}\label{lemma:mtp}
  Let $\Cn$ be $\mu$-a.s. finite. Let $\wh{X}_t$ be RRW on $H(\w)$, with initial
  distribution $\P[X_0=\cdot]$ given either by the uniform distribution
  $\Pgl[X_0=\dot]$ on $\Cn$, or by the atom at $o$, denoted by
  $\P^o[\cdot]$. Then, the annealed return probability is given by \be
  \E_\mu\left[ \P^o[\wh{X}_n = o] \right]\; =\; \E_\mu\left[
    \Pgl[\wh{X}_n=\wh{X}_0] \right].  \ee
\end{lemma}

{\em Proof:} Let $A_{vw}=P^v[\wh{X}_1=w]$ be the transition matrix of $\wh{X}$
on $\Cn\ni v,w$.  For the annealed return probability of the
vertex $o$ we have
\be
\Em \P^o\left[\hat{X}_n=o\right] \;=\; \Em A^n_{oo}
\;=\; \suli{k\in V}{}\Em \left[A^n_{oo} \frac{\chi_{(k\in \Cn)}}{\Nn}\right].
\ee

Define 
\be
f:V\times V &\to& \R,\\\la v, w\ra &\mapsto&  \E_\mu \left[A^n_{vv}
\frac{\chi_{\{w\in\mathcal{C}_v\}}}{|\mathcal{C}_v|}\right].
\ee

Applying the mass transport principle, we obtain by interchanging the arguments
of $f$,

\be
\suli{k\in V}{}\Em \left[A^n_{00} \frac{\chi_{(k\in \Cn)}}{\Nn}\right]
&=& \suli{k\in V}{}\Em \left[A^n_{kk} \frac{\chi_{(o \in
      \mathcal{C}_k)}}{|\mathcal{C}_k |}\right].
\ee

But $o\in\mathcal{C}_k \;\eqiv\; \mathcal{C}_o = \mathcal{C}_k \;\eqiv\;
\mathcal{C}_o \ni k$, so \be \suli{k\in V}{}\Em \left[A^n_{kk} \frac{\chi_{(o
      \in \mathcal{C}_k)}}{|\mathcal{C}_k |}\right] = \suli{k\in V}{}\Em
\left[A^n_{kk} \frac{\chi_{(k \in \mathcal{C}_o)}}{|\mathcal{C}_o |}\right] =
\Em \left[ \frac{1}{\Nn}\suli{k\in \Cn}{} A^n_{kk}\right]. \ee

Since $|\Cn|<\infty$, the expression on the right is $\Pgl[\wh{X}_n =
\wh{X}_0]$. These results equally apply to the continuous time RRW. In this
case, the transition kernel $A_{v,w}, (v,w \in \Cn)$ is replaced by $(\exp(\id -
A))_{v,w}$ (see [\ref{laurent}], chapter 1.3). \hfill \qed\\

This result gives the opportunity to apply the result of the former section to
the annealed return probability of RRW back to the given vertex $o$ in case of
specific random partial graphs. Invariant percolation on the Euclidean lattice
as well as on homogeneous trees in the subcritical regime shows to have a
cluster size $\Nn$ with exponentially decaying distribution [\ref{aizbar}] (see
the generalisations of this to quasi-transitive graphs in [\ref{ivto}]).\\

In order to study the return probability of RRW in this case, we consider an
invariant percolation on a transitive, unimodular graph with distribution of
the cluster size on $\{1, 2, 3, ...\}$ given by $\mu[\,\Nn = m]\,=\,C
\exp(-m/\Nh)$, for some $\Nh\in\N$. Then $\Nh\;\le\;\E[\Nn]\;\le\;\Nh+1$.\\

Let $a :=1\;+\;\frac{\sqrt{\pi\d}}{2^{1-\nu}}$, and
$r=\nu/(1+\nu)=1/((\d-1)\ln16 +1)$, using the definitions in (\ref{eq:consts}).  

Define $c:=a\left(1+e\cdot(1+(\Nh+1)^{    \frac{1}{1+2\nu}})\right)$, and
$\rho:=4^{\frac{(1+\nu)^2}{2\nu}}$.   

\begin{theo}\label{theo:app}
\;\;With the definitions given above, for  $t\ge\frac{\d}{4\hat{N}}$, it holds
\be
\E_\mu \P^o[\;\wh{X}_t\;=\;o\;]\;\le \;\P[\frac{1}{|\Cn|}]\;+\; 
  c\;t^{-\frac{1}{6}(1 \;+\; r)}\exp\{-\left( \frac{4 t}{\d
      \Nh^2}\right)^{\frac{1}{3}}\}\;+\;\sqrt{\rho}\, t^{\frac{1}{4}}e^{-
    \frac{4 \sqrt{t}}{\d \rho }}\;+\;e^{-\frac{t}{16 \wh{N}}}.
  \ee
\end{theo}

{\bf Remark:} For large values of the maximal degree $\d$, the bound retains
a prefactor $\sim t^{-1/6}$. For comparison, the trivial bound (\ref{eq:trivial})
yields a constant prefactor, along with the same exponential. Moreover, since
$\nu \sim 1/(\d-1)$, it improves with decreasing $\d$. For the constant $c$ it
holds: $c\le 8e\sqrt{d}(\Nh+1)^{\frac{1}{1+2r}}<24\sqrt{d}\hat{N}$.\\

{\em Proof:} Let $\wh{C}:=\sum_{m=1}^{\infty}
\exp(-m/\Nh)=(\exp(1/\Nh)-1)^{-1}$ be the normalising constant in
$\E[f(\Nn)]={\wh{C}}^{-1}\sum_{m=1}^{\infty}f(m)\exp(-m/\Nh)$. Write $N_o=\Nn$.
As in (\ref{eq:consts}), we let $\check{t}(N):=4^{1+\nu} N^{1-\nu}$, and
$\hat{t}(N)=4^{-(1+\nu)^2/\nu}N^4$. We note that $\check{t}(N_o) \le t$ if and only
if $N_o\le\check{N}_t $, where $\check{N}_t=4^{-\frac{1+\nu}{1-\nu}}
t^{1/(1-\nu)}$. Likewise, $\hat{t}(N_o) \le t$ if and only if
$N_o\le\hat{N}_t$, with $\hat{N}_t=\sqrt{\rho}t^{\frac{1}{4}}$. Let
$b=\frac{1-\nu}{1+\nu}=1/(1+2r) $.  By Lemma
\ref{lemma:mtp}, the expected return probability is equal to 
the expected average return probability.  By Theorem \ref{theo:haupt}, with
$I_t(\w): = \P^{\w}_{-}[\X_t=\X_0]\;-\;1/\Cn(\w)$, \;\;and
$\check{t}=\check{t}(N)$, $\hat{t}=\hat{t}(N)$ 
\be 
\E[I_t] &=&  \E[I_t| t < \check{t} \,]\P[t < \check{t} \,]\;+\;
\E[I_t|\check{t}\le t \le \hat{t} \,]\P[\check{t}\le t \le \hat{t}\,]\;+\;
\E[I_t |N_o<\hat{N}_t\,]\P[N_o<\hat{N}_t]\\ 
& \le & P[N_o> \check{N}_t\,]+
   \frac{a \wh{C}^{-1}}{\sqrt{t}^{1/(1+\nu)}} \suli{\hat{N}_t \le m\le\check{N}_t }{}
   m^b e^{-m/\Nh - (4t)/(\d m^2)}+\frac{1}{\wh{C}}\suli{m<
     \hat{N}_t}{}e^{-m/\Nh - (4t)/(\d \hat{N}_t^2)}. 
\ee

The last sum results from comparing all eigenvalues of $\exp(-t(\id-A))$ with
the second largest, as in (\ref{eq:trivial}), which is apt for $t>\hat{t}\sim N_o^4$,
as given in Theorem \ref{theo:haupt}. To the sum of the middle term we add all
missing terms corresponding to $m\in\N$ and divide it into two new
sums, the first with the terms up to $M$, and the second with terms $m=M+1,
M+2, M+3,\cdot\cdot\cdot$. 
This yields the following upper bound for $\E[I_t]$:
\be
e^{-\check{N}_t/\wh{N}}&+&
  a \,t^{-\eh\frac{1}{1+\nu}}\left( M^b e^{-\frac{4t}{\d M^2}}
    +e^{-\frac{M}{\Nh}}\wh{C}^{-1} 
    \suli{m=1}{\infty} (M+m)^b e^{-m/\Nh}\right)+\hat{N}_te^{-4t/(\d \hat{N}_t^2)}\\ \le\;\;
e^{-\frac{t}{16 \wh{N}}} &+&a\, t^{-\eh\frac{1}{1+\nu}} M^b\left(
 e^{-\frac{4t}{\d M^2}}\;\;+\;\;e^{-\frac{M}{\Nh}}\E[(1\;+\;\frac{\Nn}{M})^b]
\right)+\sqrt{\rho}\, t^{1/4}e^{-\frac{4\sqrt{t}}{\d \rho}}. 
\ee

Balancing the two exponentials inside the parentheses (by setting them equal)
would yield $M=(\frac{4}{\d}t\Nh)^{1/3}$, but this would work only for some $t$,
those for which $M$ is an integer.  However, by setting
$M_o=\lfloor(\frac{4}{\d}t\Nh)^{1/3}\rfloor$, we have $M - 1 \le M_o\le M$, and
therefore $\exp\left(-\frac{4t}{\d M_o^2} \right)\le \exp\left(-\frac{4t}{\d
    M^2}\right)$, as well as $\exp\left( -\frac{M_o}{\Nh}\right)\le
e^{1/\hat{N}}\exp\left(-\frac{M}{\Nh} \right)\le e\exp\left(-\frac{M}{\Nh}
\right)$.  Since $b<1$, it follows $(1+x)^b\le 1+x^b$ for $x\ge 0$. The
statement follows by Jensen's inequality ($x\mapsto x^b$ is concave) if $t$
obeyes $t\ge \frac{\d}{4\widehat{N}}$, implying $M\ge 1$. \hfill \qed\\

The corollary to Theorem \ref{theo:app} gives the corresponding bound in terms of
$\cp:=\E_\mu[\Nn]$ if the random events of edge removal of the Euclidean
lattice of dimension $d$ are independent.

\begin{corollary} For subcritical Bernoulli bond percolation on the Euclidean
  lattice in $d$ dimensions, the annealed return probability of the continuous
  time RRW with $\bar{a}=4\sqrt{d}$, $\rho=4^{\frac{(1+\nu)^2}{2\nu}}$, and
  $\bar{c}=160\sqrt{d}\cp^4$ has for $t\ge 
  \max\{\, d/(4\;\cp^2), \; 4^{1+\nu}\cp^{2(1-\nu)}\}$ the upper  bound
\be
\E[1/\Nn]\;\,+\,\;\bar{c}\,t^{-\frac{1}{6}(1+r)}e^{-\left(\frac{2
      t}{d} \cp^{-4} \right)^{\frac{1}{3}} }\;+\;\sqrt{\rho}\,
t^{\frac{1}{4}}e^{-\frac{2\sqrt{t}}{d  
    \rho}} \;+\;e^{-\frac{t}{32}\cp^{-2}}\;+\;a\,\cp^b
\,t^{-\frac{1}{2(1+\nu)}}e^{-\frac{2t}{d}\cp^{-4}}. 
\ee \label{cor:bperc}
\end{corollary}

{\em Proof:} The only difference if compared with Theorem \ref{theo:app}
consists of the cluster-size distribution having only an exponentially decaying
{\em tail}, and isn't necessarily of geometric (exponential) form in the first
terms. We use $\E[\cdot]:=\E_\mu[\cdot]$. The following bound is from
[\ref{grimmet}]. It is explicitly stated, when the specific exponential bound
`sets in': if $m\ge \cp^2=:\Mb$, \be \P[ \Nn \ge m] \le 2
\exp(-\frac{m}{2\Mb}).  \ee

For $m\ge \Mb$, this implies $\P[ \Nn = m] \le 2 \exp(-\frac{m}{2\Mb})$,
and $\E[f(\Nn)]=\suli{m=1}{\infty}\phi(m)f(m)$, with 
\be
\phi(m) \le C_1\chi_{{}_{<\Mb}}(m)\;\;+\;\;2 \exp(-\frac{m}{2\Mb})\chi_{{}_{\ge\Mb}}(m),
\ee

for some non-negative constant $C_1\le 1$. We mean to use Theorem
\ref{theo:haupt} in the calculation of $\E\Pgl[\X_t\;=\X_0]\;-\;
\E[1/\Nn]\le \E[ e^{-\frac{4t}{\d N_o^2}}\chi_{{}_{\N\setminus [\hat{N}_t, \check{N}_t]}}(N_o)] +
  a\sqrt{t}^{-\frac{1}{1+\nu}}\E[N_o^b  e^{-\frac{4t}{\d N_o^2}}
  \chi_{{}_{[\hat{N}_t, \check{N}_t]}}(N_o)] $,
where  $b=\frac{1-\nu}{1+\nu}$. First, we restrict $t\ge t_o$, where
$t_o=\check{t}(L)$. Then $\check{N}_t\ge L$, which, under
$t<\check{t}(N_o)$, implies $\Nn \ge L$ and $\P[t\le \check{t}(\Nn)]\le
e^{-\check{N}_t/\wh{N}}$. Furthermore, we have  
\be
\E\left[\;\Nn^b\; e^{-\frac{4t}{\d \Nn^2}}\right] \;\;\;&\le&\;\;\; \Mb^b
e^{-\frac{4t}{\d \Mb^2}} \;\;+\;\; 2\suli{m=\Mb+1}{\infty} m^b \exp(-\frac{m}{2
  \Mb} -\frac{4 t}{\d m^2}).
\ee

From here, the proof follows exactly the lines of the one of the last theorem,
except for the normalising constant $\wh{C}^{-1}$, which is missing in front of
the sum. Using Theorem \ref{theo:app}, we obtain the upper bound of
$\E\Pgl[\X_t\;=\X_0]\;-\; \E[1/\Nn]$ 
\be 
e^{-\check{N}_t/\wh{N}}\;+\;\hat{N}_t e^{-\frac{4t}{\d \hat{N}_t^2}}
&+&\frac{a}{\sqrt{t}^{1/(1+\nu)}} 
\left(\Mb^b e^{-\frac{4t}{\d \Mb^2}} \,+\,
2K^b e^{-\frac{4 t}{\d K^2}} \suli{m=1}{K} e^{-\frac{m}{2 \Mb}}\,+\,
2\suli{m>K}{} m^b\,e^{-\frac{m}{2 \Mb}}\right)\nn\\ 
\le\; e^{-\frac{\check{N}_t}{\wh{N}}} \;+\;\hat{N}_t e^{-\frac{4t}{\d
    \hat{N}_t^2}}\;&+&\;\frac{a}{\sqrt{t}^{1/(1+\nu)}} 
\left(\Mb^b e^{-\frac{4t}{\d \Mb^2}} \;+\; 4\Mb K^b e^{-\frac{4 t}{\d
      K^2}}\;+\;4\Mb K^b(1+\frac{\cp^b}{K^b})  e^{-\frac{K}{2\Mb}}\right). \ee  

Here, it was used that $\sum_{K+1}^\infty (K+m)^b\exp(-m/(2\Mb)) =
\exp(-K/(2\Mb)) K^b\E[(1\;+\;\Nn/K)^b]$ $/(\exp(1/(2\Mb)-1) $, and
$1/(\exp(1/(2\Mb)-1)\le 2\Mb$  together with Jensen's inequality and the
fact that $x\mapsto (1+x)^b$ is concave. We determine $K$ by
equating the last two exponentials:
\be 
\frac{4t}{\d  K^2}\;=\;\frac{K}{2\Mb}\;\;\;\;\eqiv\;\;\;\; K=\left(\frac{8t}{\d 
    \Mb^2}\right)^{\frac{1}{3}}.  
\ee

Again, to respect $K\in\N$, we use $\lfloor K \rfloor$, and $K-1\le \lfloor
K\rfloor \le K$.  We restrict the time parameter by
$t\ge\d/(8\Mb)$, such that $K\ge 1$. With $\d=2d$, this yields the upper bound 
\be
e^{-\frac{\check{N}_t}{\wh{N}}}\;+\;m_te^{-\frac{2t}{d m_t^2}} \;\;+\;\;  
a\,t^{-\frac{1}{2(1+\nu)}}\left( \Mb^b e^{-\frac{2t}{d\Mb^2}} \;\;+\;\;
\frac{4 \Mb^{1+b/3}}{(2d)^{b/3}} t^{b/3} (1+\sqrt{e}(1+\cp^b))
e^{-\left(\frac{2t}{d \Mb^2}\right)^\frac{1}{3}}\right).
\ee
With the definitions of $\check{N}_t$, and $\hat{N}_t$ from the proof of the
last theorem, $a<4\sqrt{d}=:\bar{a}$,\; $b=\frac{1}{1+2\nu}$, and\; 
$a(2\cp)^{2(1+b/3)}\frac{\,1\;+\;\sqrt{e}(1\;+\;   \cp^b)}{(2d)^{b/3}}\,<\,
160\sqrt{d}\cp^4=:\bar{c}$ \;\;the result follows.  \hfill \qed 

\newpage 

Now, we turn to critical Bernoulli bond-percolation in two dimensions. Here, it
is known that the connected component containing the origin is almost surely
finite [\ref{harris}].  Due to the polynomial asymptotic behaviour of the
annealed return probability [\ref{mathieu}, \ref{barlow}], the improvement
gained by the interlacing method is more significant than in the sub-critical
case.  Let $\Theta>0$ be the exponent in a scaling inequality for critical 2-d
Bernoulli bond percolation, i.e. $\P[\Nn \;\ge \; m] \le \bar{B}\cdot
m^{-1/\Theta}$ holds for some $\bar{B}>0$
([\ref{kesten}]: $\Theta$ is usually `$\d$'). 

\begin{theo}\label{theo:2dcritical} For subcritical Bernoulli bond percolation on
  the Euclidean 
  lattice in $d=2$ dimensions there is some $\Theta\ge 5$ and some $C>0$, such that
  the annealed return probability of the continuous time RRW has an upper bound
  given by 
\be 
\E\P^o[\wh{X}_t=o]\;\;\le\;\; \E[\frac{1}{\Nn}]\;\,+\,\;\frac{C}{t^{w/\Theta}},\;\;\;\;\textrm{ with }
  w \;>\; \eh\cdot(1\;+\;\frac{1}{1+9\log 16}).
\ee
\end{theo}

{\em Proof:}  We have
with $h:\R_+\to \R_+$ some given function and with 
$I_t(\w):=\P^o[\wh{X}_t=o]-1/\Nn$ and $N_o:=\Nn$
\be
\E[I_t] &=&\E[I_t|N_o \le h(t)]\P[N_o \le h(t)]\;+\;\E[I_t|N_o > h(t)]\P[N_o > h(t)]\\
&\le& \,\;\;\;\;\;\;\;\;\;\;\;\E[I_t|N_o \le
h(t)]\;\;\;\;\;\;\;\;\;+\;\;\;\;\;\;\;\;\P[N_o > h(t)], 
\ee

since $I_t(\w)\le 1$, for all $\w\in\Omega=2^E$. By Theorem \ref{theo:haupt}
the first term can be bounded by
$a\cdot\E[\left(N^{1-\nu}/\sqrt{t}\right)^\frac{1}{1+\nu} e^{-4 t/(\d
  N^2)}]$, where $a=1+\frac{\sqrt{\pi\d}}{2^{1-\nu}}$. Together with the
scaling inequality, $b=\frac{1-\nu}{1+\nu}$, and the degree $\d=2d=4$ this 
yields ($\nu=\frac{1}{3\ln 16 }\sim 0.12$)
\be
\E[I_t]\;\;\;\;\;\;\le\,\;\;\;\;\;
(1+2^\nu\sqrt{\pi}) \sqrt{t}^{-\frac{1}{1+\nu}}\;h(t)^b\;e^{-\frac{t}{h(t)^2}}
\;\;\;\;\;+\;\;\;\;\; \frac{\bar{B}}{h(t)^{1/\Theta}}.
\ee

Using $\exp(-x)\le 1/x$ for $x=t/h(t)^2$, we get by choosing $h(t)=t^w$ for
some $w>0$
\be
\E[I_t] \le
t^{\frac{w(1-\nu)-1/2}{1+\nu}\;-\;(1-2w)}\;\;\;+\;\;\;\bar{B} t^{-w/\Theta}.
\ee

Setting the exponents equal gives $w=\eh
\frac{3+2\nu}{3+\nu+\frac{1+\nu}{\Theta}}>\eh(1+ \frac{\nu}{3+\nu})$, and
$C=1+2^\nu\sqrt{\pi}+\bar{B}$. \hfill \qed\\  

{\bf Remark:} Using $I_t(\w)\le e^{-t/N_o^2(\w)}$ in the proof instead of
(\ref{eq:haupt}) produces only $w=\eh$. 

 To display the usefulness of
Lemma \ref{lemma:mtp}, we conclude this section with a lower bound for $\E[I_t]$
 for critical percolation on the binary tree, where
also a.s. $\Nn<\infty$ (see [\ref{blps2}]).

\begin{theo}\label{theo:treecritical} The annealed return probability of RRW on
  critical Bernoulli bond 
  percolation on the binary tree fulfils \;\;$\E[\P^o[\wh{X}_t=o]]\ge\E[1/\Nn]\;\;+\;\;
  c/t^{3/2}$, for some $c\ge e^{-4}/15$.
\end{theo}

{\em Proof}: By the easy part of Cheeger's inequality (see for
example [\ref{laurent}], Lemma 3.3.7), there is an upper
bound for the spectral gap $\lambda$  of
the continuous time RRW: $\lambda\le I:=\min_{A\subset \Cn: |A|\le \eh\Nn}
1/|A|$. By Lemma \ref{lemma:divide_trees} , $I\le 4(\d-1)/\Nn=12/\Nn$. By the
scaling-inequality $\phi(m):=\P[\Nn = m]\ge  
\frac{1}{10} m^{-\eh}$ (see [\ref{grimmet}], sec.10) with $I_t(\w):=\P^o[\wh{X}_t=o]-1/\Nn$
\be
 \E\left[I_t\right]=
\E\left[\frac{\Sp e^{-t\wh{\L}}}{|\Cn|}-\frac{1}{|\Cn|}\right]
\ge\E\left[\frac{e^{-t\lambda}}{|\Cn|} \right]
=\suli{m=1}{\infty}\phi(m)\frac{e^{-t \frac{12}{m}}}{m} 
\ge\suli{m\ge t}{}\frac{0.1}{\sqrt{m^3}} \frac{e^{-t \frac{12}{m}}}{m} 
\ge \frac{ e^{-12}}{15\sqrt{t}^{\,3}}.\qed
\ee

{\bf Remark:} Compare this with the asymptotic type of the RRW on the
homogeneous binary tree: $\P^o[\wh{X}_n=o]=\P^o[X_n=o]\sim \bar{\rho}^{\,n}
n^{-3/2}$, where $\bar{\rho}=2\sqrt{2}/3$ [\ref{woess}].  See also
[\ref{mathieu}], sec.4.

\newpage

\subsection{Discussion}

The random field on the edges of the graph is required to be stationary with
respect to the action of a transitive, unimodular subgroup of the group of
automorphisms.  Bernoulli percolation on $\Z^d$ is the simplest example. The
asymptotic type (as defined in [\ref{woess}], chap.  14) of symmetric random
walks with finite range on groups with polynomial growth of degree $d$ is
$n^{-d/2}$. The first who have proven accurate bounds on the return probability
of simple random walk on percolation clusters were Mathieu and Remy
[\ref{matrem}]. These were improved by Barlow [\ref{barlow}]: The return
probability of the simple random walk on the infinite percolation cluster in
$\Z^d$ has the same asymptotic quenched and annealed estimates as the random
walk on the original graph (see also Heicklen and Hoffmann[\ref{hh}], in which
the upper bound has an extra logarithmic factor). The mixing time of simple
random walk on the subgraph of an infinite percolation cluster induced by a
large box has been estimated in [\ref{matrem}] and [\ref{benmos}], the latter
of which is valid only for values of $p$ close to one. Most relevant for the
present work is the paper by Font\`es and Mathieu[\ref{mathieu}] who consider
independent random conductances of networks on $\Z^d$. They determine the
asymptotic type of quenched and annealed probabilities of a random walker.
Moreover, monotonicity of the annealed return probability under removal of
edges is proven.  Results of this type are given in [\ref{aldlyo}] for the much
more general situation of random subgraphs of unimodular transitive graphs.  \\

The present work focuses on a different subject: random walks on {\em finite}
clusters of a percolative graph. The convergence rate of random walks on finite
sets has been extensively studied. However, the finite percolative clusters
don't exhibit enough symmetry to allow the use of a general, nontrivial
estimate for the spectral gap. In general, there are no typical isoperimetric
properties, as opposed to the infinite component in the case of Bernoulli
percolation on $\Z^d$ [\ref{matrem}]. For the random walk restricted to a
finite connected component (in the subcritical case of percolation models, or
if possible conditioned on finiteness in the supercritical phase), the
 the quenched estimate is almost surely exponential.  This is not the case for the annealed return
probability.  Unlike the case of the infinite cluster, there are arbitrarily
large clusters with bounded edge-connectivity related to a small spectral gap.
By the annealing, these large clusters contribute to some degree, and it is
natural to ask whether this influences the expected return probability to the
extent of a retarding effect. A result by I.Benjamini and O.Schramm used in
[\ref{hh}] (Theorem 3.1) on finite subgraphs of percolation clusters implies
that under removal of edges the expected return probability is non-decreasing.
In particular, it is of interest by how much at most the clusters with a
particularly large relaxation time (in terms of the cluster size), such as
paths, slow down the walk, i.e. decrease the asymptotic decay of the return
probability. This question is answered by the more difficult upper bound
(\ref{theo:partial}) which results from the second of the two inequalities in
(\ref{eq:il_ineqs}). Furthermore, it is of interest if the quality of
stationarity in an invariant percolation leads to a speed up of the convergence
of the simple bounds resulting from assuming for each connected component the
geometry with the longest relaxation time. This is answered positively in the
case of unimodular graphs, allowing for the application of Lemma
\ref{lemma:mtp}, which states the equivalence of the expected return
probability with the expected average return
probability. For this quantity Theorem \ref{theo:app} testifies an improved polynomial
prefactor $\sim t^{-\frac{1}{6}(1+\nu/(1+\nu))}$ (see (\ref{eq:consts}) for the
definition of $\nu$).\\  

Note that the return probability of the simple random walk on (deterministic)
Cayley graphs of polycyclic groups with exponential growth (see [\ref{woess}],
chapter III., 15), lamplighter groups [\ref{revelle}], and Diestel-Leader
graphs [\ref{barwoe}] without drift also show the characteristic $e^{-c
  n^{1/3}}$ for some $c>0$ .  Our estimate concerns the finite percolative
subgraphs of all $\d$-regular graphs which have a {\em unimodular}, transitive
subgroup $\Gamma$ of the automorphism group.  Since the lamplighter group is
amenable, it is unimodular [\ref{sw}], while among the Diestel-Leader graphs
are examples of transitive graphs with a non-unimodular automorphism group
(e.g.  [\ref{barwoe}]) and they don't fall into the present scope. \\

In [\ref{kirsch}], Lemmata 2.7, 2.9, it is shown, that the integrated density of
states IDS (cumulative spectral function) $E\mapsto N_N(E)=\E(e_o, \id_{[0,
  \infty)}(E-L) e_o)$ (in the notation of [\ref{kirsch}]) of combinatorial
Laplacians $L$ (the subscript $N$ stands for `Neumann Laplacians'; see the remark
after Theorem \ref{theo:partial}) of subcritical Bernoulli percolation
graphs on a Euclidean lattice $\Z^d$ has `Lifshits tails' of a specific form:
there are constants $\a^-, \a^+$, such that 
\bel 
\exp(-\a^- E^{-{1/2}})\;\;\le \;\;N_N(E) \;-\;N_N(0) \;\;\le\;\;\exp(-\a^+
E^{-{1/2}}), \label{eq:alpha}
\eel 

where $C_o$ is the finite percolation cluster containing the origin of the
lattice.  No further specification of these constants is given. For the upper
bound, we specify a value of $\a^+$, for which this estimate is valid.\\

As was stated in section \ref{sec:interlacing}, the combinatorial Laplacian $L$
of a graph remains invariant under decoration with loops, so $L=\wh{L}\equiv
\d\id - \wh{\mathcal{A}}$, where the loops are those of the regularisation of
definition 1.1. The transition matrix of the RRW on the graph is $A=\d^{-1}
\wh{\mathcal{A}}$, and the normalised Laplacian $\wh{\mathcal{L}}$, of the
regularised graph is given by $\wh{\mathcal{L}}=\id - \d^{-1}(\d\id - \wh{L})=
\d^{-1}L$. Therefore, spectral estimates about the RRW on a graph always also
apply to the scaled Laplacian spectrum of the graph, as in Lemma
\ref{lemma:li}. In case of Bernoulli percolation on $\Z^d$, the scaling factor is
$\d=2d$.\\ 

The connection with the integrated density of states is given by the following
well-known lemma, which states, that the IDS is the expected continuous time
return probability of the RRW on $C_o$, the percolation cluster containing the
origin.

\begin{lemma} Let $\wh{X}_t$ be RRW on the finite, subcritical percolation
  cluster $\Cn$ of Bernoulli bond percolation on $\Z^d$ containing the origin. Then
\be 
\ili{0}{\infty} e^{- t E} dN_N(E)\;\;=\;\;
\E_{\mu}\P^o[\wh{X}_{2d t}=o]\;\;-\;\;\E_{\mu}[\frac{1}{|C_o|}]. 
\ee \label{lemma:lap}
\end{lemma}

{\em Proof:} Due to ergodicity ([\ref{kirsch}], Lemma 1.12), we have
weak convergence of  
\be 
N_N^n(E):=\E_\mu \left[\frac{1}{|\Lambda_n|}\Sp \id_{[0,\infty)}(E-L_n) \right]
\ee 
to $N_N(E)$ as $n\to\infty$, where $L_n$ is the Laplacian belonging to the
subgraph of the subcritical percolation graph induced by (restricted to) $
\Lambda_n\;:=\;\{\;-n\;+1\;\;, \;...\;,\;\;n\;\}^d$. Since the function $x\mapsto 
\exp(-x)$ is bounded and continuous on $\R_+$, this implies convergence of
$\int_0^\infty e^{- t E} dN^n_N(E)$ to $\int_0^\infty e^{- t E} dN_N(E)$,
as $n\to\infty$.  Moreover, since there are only finitely many
configurations on the edges of the finite graph, the expected value $\E_\mu$ is
just the arithmetic mean over these configurations, so:
\bel \ili{0}{\infty} e^{- t E} dN^n_N(E) = \E_\mu\left[ \ili{0,}{4d}e^{-E t}
  \frac{1}{|\Lambda_n|}   \,d\,\Sp(\id_{[0,\infty)}(E-L_n))\right].  
\label{eq:fin_lap} \eel

The integration bound follows from $\sigma(L^n)\subset[0,4d]$, see
[\ref{kirsch}], remark 1.10.(ii).  Let $\Pn$ be the $l^2(\Z^d)$-projector onto
$l^2(\Lambda_n)$. Since $L_n=\Pi_{\Lambda_n}L\Pi_{\Lambda_n}$
is of finite rank, by the spectral theorem for real, symmetric matrices, the right
hand side equals $\E_\mu[ (1/|\Lambda_n|)\Sp\exp(- tL_n )]$, the average return
probability of the RRW on the subgraph of the percolation graph induced
by $\Lambda_n$.\\

We now use Lemma \ref{lemma:li}.  Considering the pertubation $S_n$ in $L_n=\Pn L \Pn\;+\;S_n$,
it is clear that $S_n \ge 0$ and rank$(S_n)\le |\partial \Ln|$ (the edge
boundary of $\Lambda_n$). Then, by the second inequality of (\ref{eq:li}),
\be
\E_\mu (e_o, \exp(-tL) e_o)\;\;=\;\;
\E_\mu \frac{1}{|\Lambda_n|}\Sp\Pi_{\Lambda_n}\exp(-tL)
\;\;\le\;\;\E_\mu \frac{1}{|\Lambda_n|}\Sp\Pi_{\Lambda_n}\exp(-tL_n).
\ee

Since $e^{-x}\;-e^{-y}\;\le\; y\;-\;x$ for $0\le y\le x$, the right hand side
can be further approximated:
\be
\E_\mu
\frac{1}{|\Lambda_n|}\Sp\Pi_{\Lambda_n}\left(\exp(-tL_n)\;-\;\exp(-tL)\right)
\;&+&\;\E_\mu\frac{1}{|\Lambda_n|}\Sp\Pi_{\Lambda_n} \exp(-tL)\\
\;\;&\le&\;\;t\frac{|\partial \Lambda_n|}{|\Ln|}\;\;+\;\;\E_\mu (e_o, \exp(-tL) e_o).
\ee

Therefore, due to the amenability of the Euclidean lattice, 
\bel
\lim\limits_{n\to \infty} \E\frac{1}{|\Ln|}\Sp \exp(-tL_n)\;\;=\;\;\E(e_o,
\exp(-tL) e_o).\label{eq:L_n<->L}
\eel

With the arguments given above, this implies $
\int\limits_0^\infty e^{-tE} dN_N(E)\;\;=\;\;\E_\mu(e_o, e^{-tL} e_o)$.

Now,
$A=\id-\wh{\mathcal{L}}=\d^{-1} \mathcal{\wh{A}}$ is the transition operator of
the RRW on the random connected components. With $L=\widehat{L}=\d
\wh{\mathcal{L}} = \d (\id - A)$, the right hand side of (\ref{eq:L_n<->L}) is
$\E_{\mu}\P^o[\wh{X}_{\d t}\;=\;o]$, the expected return probability
of the RRW on $\Cn$, where $\d=2d$.  By remark 1.15, (ii), $N_N(0)=\E[1/|\Cn|]$
(see also [\ref{grim_num}]). \hfill \qed\\

{\bf Remarks:} Due to results in [\ref{lemuve}] (see theorem 2), the
convergence of $N_N^n$ to $N_N(\cdot)$ is even uniform. The application of
Lemma \ref{lemma:li} is similar to the application of Theorem 3.1 to the proof
of Theorem 8.1 in [\ref{hh}].\\

Now, Theorem \ref{theo:app} will be applied in such a way as to show: $\alpha^+\ge
\frac{4}{3\sqrt{3}}\cp^{-4} $:

\begin{theo} The integrated density of states, $E\mapsto N_N(E)$, of Laplacians
  belonging to subcritical Bernoulli percolation on the Euclidean lattice in
  $d$ dimensions fulfils with $0<E\le \wh{E}$, s.t.
  $\wh{E}=\min\{ \;\frac{d^{1/3}\cp^{4/3}}{3\cdot 2^{1/3}}, \;
  \frac{\sqrt{\rho}/6}{d^{2/3} \cp^{4/3}}, \frac{2^{1/3} \rho^{4/3}}{768 \cp^4},\;
  \frac{8/3}{(\d\rho\cp)^{4/3}}\}$, $\bar{c}=160\sqrt{d}\cp^4$, $\rho=4^{\frac{(1+\nu)^2}{2\nu}}$
\be N_N(E) \;\;-\;\;N_N(0)\;\;\;\le\;\;\;A E^{} \exp(-\a E^{-1/2})\;\;+\;\; B
  E^{-3/8}\exp(-\beta E^{-3/4}),
\ee

where $\a = \frac{4}{3\sqrt{3}}\cp^{-4}$, $A= 3\sqrt{3}\,
2^{-\frac{1}{6}(1+r)} \bar{c} \left(d \,\cp^4\right)^{\frac{1}{6}(1+r)}$, 
$B=3\sqrt{\rho}\frac{(2/d)^{1/4}}{3^{3/8}}\cp^{-1}$,
$\beta=\frac{2\sqrt{2}}{(\sqrt{3}d)^{3/2}}\frac{1}{\rho \cp^{2}}$.\\ 
\end{theo}

{\em Proof:} We begin with a lower bound of the Laplace transform of the IDS:
\be
\ili{0}{\infty}e^{-t E}dN_N(E)\;
\;\;\ge\;\;\ili{0}{\e}e^{-t E}dN_N(E)\;\ge\;\;e^{-t \e}(N_N(\e) \;-\;N_N(0)).
\ee

Due to Lemma \ref{lemma:lap}, we may apply Corollary \ref{cor:bperc} to obtain
an upper bound. From the corollary, we see that  for all $t\ge t_0:=2/(d\cp^4)$
with $r=\nu/(1+\nu)$
\be
N_N(\e)-N_N(0)&\le&e^{t\e}\left( \bar{c}t^{-\frac{1}{6}(1+r)}e^{-\left(\frac{2
        t}{d} \cp^{-4} \right)^{\frac{1}{3}} } + \sqrt{\rho}
  t^{1/4}e^{-\frac{2\sqrt{t}}{d \rho}}+e^{-\frac{t}{32}\cp^{-2}}+a\cp^b
  t^{-\frac{1}{2(1+\nu)}}e^{-\frac{2t}{d}\cp^{-4}} \right)\;\;\;\;\;
\ee

It is easy to check, that for $t\ge t_1:=\max\{ 4\sqrt{d/\rho}, \;\cp^8/\rho,\;d^2\rho^2/16
\}$, the third and fourth term in the parentheses is smaller or equal to the
second, yielding the factor 3 in constant $B$. The first term is
$\bar{c}t^{-\frac{1}{6}(1+r)}\exp\left(\e t\;\;-\;\;\bar{\a}^{\frac{1}{3}}
  t^{\frac{1}{3}} \right)$, where $\bar{\a}=\frac{2}{d 
  \,\cp^{4}} $. We set $c\, \e\, t = \bar{\a}^{1/3} t^{1/3}$, for some
$c>1$. Then,  $t =  \frac{\sqrt{\bar{\a}}}{ (c\, \e)^{3/2}}$. This yields for the exponent
\be
t\,\e - c\,\e\,t = -\frac{c\,-\,1}{c^{3/2}}\frac{\sqrt{\a}}{\sqrt{\e}}.
\ee

We optimise the result by setting $c=3$. The upper bound of $N_N(\e)-N_N(0)$ follows
by inserting the resulting $t$ into the inequality. Restricting
$t=t(\e)\ge \max\{ t_0, \;t_1\}$ gives the range $[0,\wh{E}]$ for
$\e$, in which the first two terms of the right hand side dominate.  \hfill \qed\\  

{\bf Remark:} A different situation prevails in the supercritical case. As
shown in [\ref{muesto}], the asymptotics of the SRW on the infinite cluster
derived in [\ref{barlow}] correspond with and characterise a different
behaviour of the spectrum of the graph Laplacian in this case. The exponential
decay of subcritical Bernoulli percolation clusters on quasi-transitive
graphs has been shown in [\ref{ivto}].

\section{Interlacing as a new technique for comparison theorems}  

We now introduce a general method to derive bounds for the average return
probability of the RRW on a finite graph $H$. It consists of four steps: i.)
Interpreting $H$ as a connected component of a larger graph, by adjoining an
additional connected component ($\tilde{H}$); 

ii.) Comparing the graph $H+\tilde{H}$ with a transformation thereof, called
$H'$, on which the RRW has known spectrum (here, interlacing is used);
\;\;\;iii.)  Performing an estimation of the RRW on $H'$; \;\;\;iv.) Optimising
the result, e.g.  by choosing the optimal size of $\tilde{H}$.

\subsection{Proof of the main result}

We begin by noting that for RRW $\wh{X}_\cdot$ on a finite, simple, connected
graph $H$, 

\be
\Pgl[\X_t=\X_0]&=& \frac{1}{N}\Sp[\exp -t(1-A)],
\ee

which is the return probability of the RRW on $H$ in continuous time (Observe:
$A=\d^{-1}\wh{\mathcal{A}}$, where $\wh{\mathcal{A}}$ is the adjacency matrix of
  the regularized graph).  We then have

\be
\Pgl[\X_t=\X_0]\;-\;\frac{1}{N}\;\;&=& \;\;\frac{1}{N}\suli{j=2}{N}{\a}_j^t,
\ee
where $\a_j=\exp(-(1-\beta_j))$, with $\beta_j$ the $j$-th element of the
spectrum of $A$, ordered in a non-increasing fashion: $1=\beta_1 > \beta_2 \ge
\beta_3\ge \cdot\cdot\cdot \ge \beta_{N}$. Furthermore, $\a_j^t\equiv
(\a_t)^t$, and $1=\alpha_1 > \alpha_2 \ge \cdot\cdot\cdot \ge \alpha_N>0$. \\ 

The idea is to compare RRW on $H$ with RRW on another graph, a reference graph,
about which the spectrum of the transition kernel is known. Interlacing can be
usefully employed, whenever the rank of the transformation of $H$ into the
reference graph is small compared to the order of $H$. We will use the finite
path $P_{N'}$ (for some size $N'$) as a reference graph. In general, the
transformation of $H$ into $P_{N'}$ will entail rearranging a number of edges
comparable to $N$, the order of $H$. This implies that the rank $K$ of the
transformation of the corresponding transition kernels does {\em not} fulfil
$K\ll N$.\\ 

To still make interlacing applicable, we remedy this circumstance by adjoining
another graph to $H$ as a different connected component, namely, the finite
path $\tilde{H}=P_{\tilde{N}}$. We then consider the transformation of the
reducible graph $H\;+\;P_{\tilde{N}}$ into $H':=P_{N'}$ as the relevant
transformation $T$ for estimating the eigenvalues $\beta_j$ by interlacing. The
rank of $T$ divided by the total number of vertices, $N'=N+\tilde{N}$ will be
arbitrarily small if $\tilde{N}$ is chosen large enough.\\

\begin{figure}[htb]
\centerline{
\epsfxsize=10cm
\epsffile{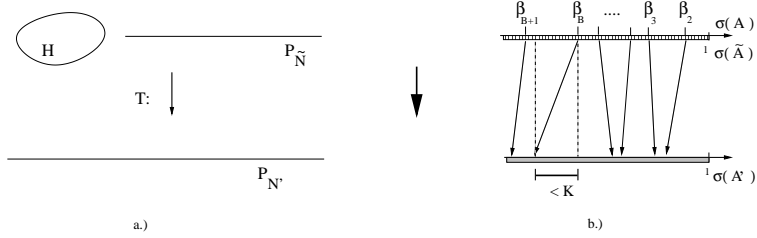} \hskip 1.5cm}
\caption{a.) Transformation of a disconnected graph $H + P_{\tilde{N'}}$ into the
  path $P_{N'}$; \;\;b.) schematic of the corresponding change of the spectrum:
  the imbedded $\beta_j=\bar{\beta}_{\iota(j)}$ shift across less than $K$
  intervals of eigenvalues of RRW on $P_{N'}$ ((\ref{eq:il_ineqs}) in Theorem
  \ref{theo:partial}).}\label{fig:trans}
\end{figure}

In other words, we consider $H$ as part of a larger graph by adjoining
$P_{\tilde{N}}$, which will not be changed during the transformation (see Fig.
\ref{fig:trans}, a.). The size $\tilde{N}$ of the adjoined component will then
be taken large enough in order to guarantee that the number of edges which have
to be rearranged to obtain $P_{N'}$ is small in comparison to $N'$. This
enables to control the change of the spectrum of the transition kernel of RRW
on $H$ by interlacing, if considered as a perturbation of RRW on $P_{N'}$.
Enlarging $\tilde{N}$ has the effect that $N'=N+\tilde{N}$ becomes larger, and
the intervals between successive eigenvalues of RRW on $P_{N'}$ smaller. \\

In the case $N=1$, there is nothing to prove, so  assume $N\ge 2$. We
proceed in four steps: \\ 

i.) To abbreviate notation, let $\ba_i=\exp(-(1-\bar{\beta}_i))$, and
$\bar{\beta}_i$ be the $i$th eigenvalue of $A\bigoplus \tilde{A}$ with
$\tilde{A}$ the transition kernel of a RRW on the finite path $P_{\tilde{N}}$
of length $\tilde{N}$.  Taking some value $B\in\{0, ..., N-2\}$, we have
\bel
\Pgl[\X_t=\X_0]\;-\;\frac{1}{N}\;\;&\le&\;\;\frac{B}{N}\a_2^t\;\;+\;\;
\frac{1}{N}\suli{i=\iota(B+2)}{\iota(N)} \bar{\a}_i^t,\label{eq:d'abord}
\eel

where $\iota:\{1, ..., N\} \to \{1, ..., N'\}$ is the injective map assigning
$\beta_j$ its index in the spectrum of $A\bigoplus\tilde{A}$, which is the
transition kernel of the reducible RRW on $H+P_{\tilde{N}}$. We choose the
eigenvalues $\bar{\beta}_i$ also enumerated in a non-increasing way. This implies
the following properties of $j\mapsto \iota(j)$:
\bel
\iota(j) \ge j,\;\;\;\;\;\;\iota(N)\le N',\;\;\;\;\;\;j\in\{1, ..., N\}.\label{eq:iota}
\eel

ii.) We now perform the transformation $T$ of the state space, the disconnected
graph $H+P_{\tilde{N}}$, to a connected graph $G'$ of equal order
$N'=N+\tilde{N}$, namely, the finite path $P_{N'}$.\\

We may assume $H$ to be a tree. If $H$ has cycles, we take out edges, until a
spanning tree is achieved. This can only increase $\a_j$. Then, the
transformation $T$ entails pairs of operations, each including a removal and an
insertion of an edge. Note, the total number of insertions of edges doesn't
exceed $N-1$.  Correspondingly, the rank $K$ of the negative part of the
perturbation $S$ in $A'=A\bigoplus \tilde{A} + S$ is bounded by $N$ (see
Fig. \ref{fig:trans} b.). \\

The following is the crucial step involving interlacing: let
$\ap_i=\exp(-(1-\beta'_i))$, where $\beta_{i+1}'=1-\frac{2}{\d}(1-\cos(\pi
i/N'))$ is the $1+i$'th eigenvalue in $\sigma(A')$. Now, we ask $\tilde{N}$ (or
$N'$) to be sufficiently large, such that the interval $(\beta_{B+2}, 1)$
contains enough elements of the spectrum of $A'$, such that
$\bar{\alpha}_{\iota(B+2)}\le \a'_{K+X}$, for some integer $X \ge 1$, which we
specify, later. This is always possible, since $|\sigma(A')\bigcap [1-\epsilon,
1]|=O(N')$, for any $\epsilon>0$, as $N'\to\infty$. In other words, we set $N'$
appropriately, such that the condition

\bel\label{eq:claim} \beta_{B+2} \;\;\;< \;\;\;\beta'_{K+X} \eel

holds.  Then, also $\bar{\a}_{\iota(B+2)} < \a'_{K+X}$, and $\iota(B+2)\ge
K+X+1.$ So, by Theorem \ref{theo:partial}, the second term of the right hand
side in (\ref{eq:d'abord}) can be further bounded from above, and using
condition (\ref{eq:claim}), and (\ref{eq:iota}): 
\bel
\frac{1}{N}\suli{i=\iota(B+2)}{\iota(N)} \bar{\a}_i^t \;\;\le\;\;
\frac{1}{N}\suli{i=\iota(B+2)}{\iota(N)} (\ap_{i-K})^t.\;\;\le\;\;
\frac{1}{N}\suli{i=X+K+1}{\iota(N)} (\ap_{i-K})^t.\;\;\le\;\;
\frac{1}{N}\suli{i=X+1}{N'} (\ap_{i})^t.\label{eq:il} \eel

iii.) Finally, approximating the cosine by a quadratic polynomial ($\cos(\pi x) \le
1-2x^2$, for $x\in[0,1]$), and using $\int_z^\infty \exp(-x^2)dx \le
(\sqrt{\pi}/2) \exp(-z^2)$  (see [\ref{laurent}],
Example 2.1.1), this can be further bounded from above by

\be 
\frac{1}{N}\suli{i=X+1}{N'}\exp(- t\frac{4}{\d} \frac{i^2}{N'^2}) \;\;\;\le\;\;\;
\frac{N'}{N} \int\limits_{X/N'}^\infty \exp(-t\frac{4}{\d}y^2)dy \;\;\;\le\;\;\;
\frac{\sqrt{\pi \d}}{4\sqrt{t}}\frac{N'}{N} \exp(- t\frac{4}{\d} \frac{X^2}{N'^2}).
\ee

Putting it all together yields
\bel
\Pgl[\X_t=\X_0]\;-\;\frac{1}{N}\;\;&\le& \;\; \frac{B}{N}\exp(-\frac{4t}{\d N^2})\;\;+\;\;
\frac{\sqrt{\pi \d}}{4\sqrt{t}}\frac{N'}{N}\exp(-t\frac{4  X^2}{\d
  N'^2}). \label{eq:result1}
\eel

iv.) What remains is to select the parameters $B, N'$, and $X$, such that, on
the one hand, condition (\ref{eq:claim}) is fulfilled, and, on the other, the
bound becomes optimal.\\ 

First, to meet condition (\ref{eq:claim}), we use Theorem \ref{theo:higher} in
the guise of  Corollary \ref{cor:higher}. From (\ref{eq:higher}), we see that it is fulfilled, if 
$\beta^o_{[((B+1)/2)^\nu]+1} < \beta'_{K+X}$. Equivalently, since

$\beta^o_j=1-\frac{2}{\d}(1-\cos(\pi j /N))$, and
$\beta'_j=1-\frac{2}{\d}(1-\cos(\pi j /N'))$, the claim (\ref{eq:claim}) is true if
\be
\frac{[((B+1)/2)^\nu ]+1}{N} > \frac{K+X}{N'}.
\ee

Using $[x]+1> x$,  and choosing $X=N$, it shows that this condition is met if 
\be 
N' := 2^\nu \frac{N(K+N)}{(B+1)^\nu}.
\ee

Now, we use the observation made above that $K\le N$, and set $B:=[xN]$ with
\;\;$x\in [0,\eh)$ to give an upper  bound of the right hand side of
(\ref{eq:result1}), using $ xN-1\;< \;B\;\le\; xN$:

\be
\Pgl[\X_t=\X_0]\;-\;\frac{1}{N}\;\;&\le& \;\;\left(x +
  \frac{\sqrt{\pi \d}}{2^{1-\nu}\sqrt{t}}\frac{N^{1-\nu}}{x^\nu}e^{-\frac{4t}{\d
      N^{2}}\left( 2^{-2(1+\nu)}  x^{2\nu} N^{2\nu} -1\right)}\right) e^{-\frac{4t}{\d N^2}}.  
\ee

By setting $x=(N^{1-\nu}/\sqrt{t})^{1/(1+\nu)}$ with
$t>\check{t}=4^{1+\nu}N^{1-\nu}$ (making $x<\eh$), the exponential inside the
parentheses remains less than one for $t\le\hat{t}=4^{-(1+\nu)^2/\nu}N^4$, whence
the exponent remains negative. Inserting this bound, we obtain the result claimed by Theorem
\ref{theo:haupt}.\hfill \qed

\subsection{A lemma on $\d$-bounded trees }

\begin{lemma}\label{lemma:divide_trees}
  Every finite tree $T=\la V, E\ra$ of order $N=|V|$ with at least two
  vertices and largest degree $\d\ge 2$ can be divided into two subtrees $T_1$ and
  $T_2$ by removal of an edge $e\in E$, such that the ratio $m$ of the
  cardinalities of the vertex sets of $T_1$ and $T_2$ obeys: \be 1 \;\;\le\;\;
  m \;\;\le\;\; 4(\d-1) -1.  \ee 
\end{lemma}

Before giving the proof, we consider some definitions concerning any finite
tree $T=\la V, E\ra$: Let a path in $T$ of length $n\in \N$ be an $n+1$-tuple
of vertices $\bar{v}:=\la v_0, v_1, ...., v_n\ra\in V^{n+1}$ with $\{v_{k-1},
v_k\} \in E$, for $k\in \{1, ..., n\}$. Let $u\leftrightarrow_{\bar{v}} w$ for
$u,w\in V$ signify that there is a path $\bar{v}$ in $T$ of length $n$ for some
$n\in \N$, such that $v_0=u$, and $v_n=w$. Say that the path $\bar{v}$ crosses
the edge $e\in E$ if $e = \{v_{k-1}, v_k\}$, for some $k\in \{1, ...., n\}$.

\be 
V_{u,e}\;:=\;\{\;w\in V\;:\;\;u\leftrightarrow_{\bar{v}} w \textrm{ for some
  path }\bar{v} \textrm{ which does not cross } e\;\}. 
\ee

In other words, $V_{v,e}$ ist the vertex-set of the subtree containing $v$
which results from removing edge $e$. Furthermore, call an edge $e_c=\{v, w\}
\in E$ to be a {\bf central edge} of $T$ if

\bel
|V_{v, e_c}| = \min\limits_{ e\in E} \max\limits_{u\in e} \left|V_{u,
    e}\right|.   \label{eq:defec} 
\eel

The minimum always exists, as long as $V_{v, e_c}$ is non-empty. This is
the case for trees of order at least two. $e_c$ is however not unique, as the
example of a ball in a regular tree shows.  Removing one of the central edges
from a tree results in two subtrees, the larger of which is as small as
possible.  For the order of the subtree $T|V_{u,e}$,  write
$N_{u,e}:=|V_{u,e}|$.\\

{\em Proof:} (Lemma \ref{lemma:divide_trees}) In the trivial cases $N=2, 3$,
the ratio $m$ is $1, 2$, respectively. Let $N\ge 4$. Let $\{ u, v \}:=e_c$ be
a central edge of $T$ with $N_{v, e_c} \ge N/2$. Consider an incident edge
$e'=\{v, w\}$ of $e_c$ (where $e'\neq e_c$), such that $N_{w,e'}=\max\{\;
|V_{r,\{v,r\}}|\;:\;\{ v,r \} \in E \}$. Such an edge always exists, since
$T_{v, e_c}$ has more than one vertex. Then, by the Pigeon-hole principle
([\ref{doug}], Lemma 1.3.10), we have

\bel
N_{w,e'} \ge \frac{N_{v,e_c}-1}{\d -1}.\label{eq:ineqtree}
\eel

By the definition of the central edge,  $N_{u,e_c}\le N_{v, e_c}$, and
$N_{v,e'}> N_{w,e'}$. On the other hand, 

\be
\frac{N_{v, e'}}{N_{w, e'}} \;=\; \frac{N\;-\;N_{w,e'}}{N_{w,
    e'}} \;=\; \frac{N}{N_{w,e'}} - 1 \;\le\; (\d - 1)\frac{N}{N_{v,e_c}-1} - 1
\;\le\; 4(\d-1) -1.
\ee

The first inequality is (\ref{eq:ineqtree}). The second inequality follows
from the definition (\ref{eq:defec}) of $e_c$, giving $N_{v,e_c}\ge
\frac{N}{2}$, and from the assumption $N\ge 4$. \hfill \qed

\begin{figure}[htb]
\centerline{
\epsfxsize=6.5cm
\epsffile{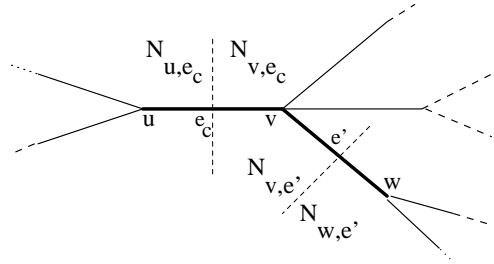} \hskip 1.5cm}
\caption{Sketch of the part of a tree $T$ around a central edge $e_c=\{u,
  v\}$ (see proof of Lemma \ref{lemma:divide_trees}), where the ratio $N_{v,e_c}/N_{u,e_c}$
  of the orders of the two subtrees $T(v,e_c)$ and $T(u,e_c)$ is the smallest possible
  above one.}\label{fig:central_e}
\end{figure}

\section{Outlook and Acknowledgement}

In a forthcoming paper, the method of obtaining comparison theorems for return probabilities
with interlacing is applied to random walks on infinite graphs.\\

My special thanks go to I. Krasovsky, who initiated the method with his work on
pertubation theory [\ref{igor}], to W. Woess, for his support with the idea of
generalizing the topic of my thesis, and to R. Lyons, P. Mathieu and W. Kirsch
for helpful remarks. The author is suppported by FWF (Austrian Science Fund),
project P18073.

\section{Bibliography}

\begin{enumerate}

\item{M. Aizenman, D. J. Barsky: `Sharpness of the Phase Transition in
    Percolation Models', CMP 108, 489-526, 1987}\label{aizbar}
 
\item{D. Aldous, R. Lyons: `Processes on unimodular random networks', 
 Electron. J. Probab.  12  (2007), no. 54, 1454--1508}\label{aldlyo} 

\item{T. Antunovi\'c, I. Veseli\'c: `Sharpness of the phase transition and
    exponential decay of the subcritical cluster size for percolation on
    quasi-transitive graphs', Journ. Stat. Phys. 130(5):983-1009, 2007  }\label{ivto} 

\item{T. Antunovi\'c, I. Veseli\'c: `Spectral asymptotics of percolation
    hamiltonians on amenable Cayley graphs', In: Proceedings of OTAMP
    2006. Operator Th.: Adv. Appl. 2007}\label{antves}  

\item{M. Barlow: `Random walks on supercritical percolation clusters',
    Ann. Probab.  32,  no. 4, 3024--3084., (2004)}\label{barlow} 

\item{I. Benjamini, R. Lyons, Y. Peres,  O. Schramm: `Group-invariant
    percolation on graphs', Geom. Funct. Anal. 9 (1999), 29--66.}\label{BLPS}

\item{I. Benjamini, R. Lyons; Y. Peres, O. Schramm:  `Critical percolation on
    any nonamenable group has no infinite clusters',  Ann. Probab.  27,
    no. 3, 1347-1356, 1999}\label{blps2} 

\item{L.Bartholdi, W.Woess: `Spectral computations on lamplighter groups and
    Diestel-Leader graphs', J. Fourier Analysis Appl., 11, 2, 175--202, 2005
  }\label{barwoe} 

\item{I.Benjamini, E.Mossel: `On the mixing time of simple random walk on the
    super critical percolation cluster',  Prob. Theory Rel. Fields  125
    (2003),  no. 3, 408--420}\label{benmos} 

\item{B. Bollob\'as, , V. Nikiforov: 
`Graphs and Hermitian matrices: eigenvalue interlacing',
Discrete Math. 289 (2004), no. 1-3, 119--127.}\label{bolobas}

\item{M. Brown: `Interlacing Eigenvalues in time reversible Markov chains ',
    Math. Oper. Res., Vol. 24, No.4, Nov. 1999}\label{bro} 

\item{D. Chen, Y. Peres, (G.Pete): `Anchored Expansion, Percolation and Speed',
  Ann. Probab. 32, no. 4, S. 2978-2995, 2004}\label{cheper2} 

\item{G. Chen, G. Davis, F. Hall, Z. Li, K. Patel, M. Stewart: `An interlacing
    result on normalized Laplacians', SIAM J. Discrete Math. , Vol. 18, No. 2,
     353-361, 2004}\label{chen_und_co}

\item{ F.R.K. Chung: `Spectral Graph theory', CBMS, Am. Soc. Math. Sc.,
    No. 92, 1997}\label{chung}

\item{D.Cvetkovi\' c, M.Doob, H.Sachs: `Spectra of graphs', Verl. d.
    Wiss.,  Berlin 1982}\label{cds}

\item{P. Diaconis, D.Stroock:  `Geometric bounds for eigenvalues of Markov
    chains', Ann. Appl. Probab. 1, no. 1, 36--61., 1991}\label{diastr}

\item{W. Feller: `An introduction to probability theory and its applications',
    Vol. I., 3ed ed. John Wiley \& Sons,  436, 1968}\label{feller}

\item{M.E.Fischer, J.W.Essam:`Some cluster size and percolation problems',
    J. Math. Phys. 2, 609, 1961}\label{fishess}

\item{L.R.G. Font\`es, P. Mathieu: `On symmetric random walks with random conductances',
Prob. Theory Rel. Fields, vol. 134, No 4, 565-602, 2006 }\label{mathieu}

\item{C. Godsil, G. Royle: `Algebraic graph theory', Springer, 2001
    chapter. 9}\label{godsil} 

\item{G.R.Grimmet: `On the number of clusters in the percolation
    model', J.London Math.Soc. (2), 13 (1976),
    346-350}\label{grim_num}

\item{G.R.Grimmet: `Percolation', chapter 4,  132, Springer, 2000}\label{grimmet}

\item{J. Guo : `Bounds on the $k$'th Eigenvalues of Trees and Forests',
    Lin. Alg. App., 149:19-34, 1991}\label{guo}

\item{J. Guo : `The $k$'th Laplacian eigenvalue of a tree', J. Graph Th.,
  54, 1, 51 -57, 2006}\label{guo2}

\item{W.H. Haemers: `Interlacing eigenvalues and graphs', Lin.A.App.
    226/228, 593-616, 1995 }\label{haemers} 

\item{T.E. Harris: `A lower bound for the critical probability in a certain
    percolation process', Proc. Cambr. Philos. Soc. 56, 13-20}\label{harris}

\item{D. Heicklen, C. Hoffmann: `Return times for simple random walks on
    percolation clusters', Electronic Journal of Prob. 10, No. 8, 250-320,
    2005}\label{hh} 

\item{J. van den Heuvel: `Hamiltonian cycles and eigenvalues of graphs', Lemma 2,
    Lin. Alg. Appl. 226/228, 723-730, 1995}\label{heuvel}

\item{R.Horn, Ch.Johnson: `Matrix Analysis', Cambridge
    Univ. Press, 1985}\label{hj} 

\item{H. Kesten: `Scaling Relations  for 2D-Percolation', CMP 109, 109-156,
    1987}\label{kesten}

\item{J.F.C. Kingman: `',  J. R. Stat. Soc. B 30,  499-510, 1968 }\label{king}

\item{W.Kirsch, P.M\"uller: `Spectral properties of the Laplacian on
    bond-percolative graphs': Math. Z. 252, 899-916, 2006}\label{kirsch}

\item{I.Krasovsky, V.I. Peresada: 'Principles of finite-dimensional Pertubation
    theory', Low.Temp. Phys. 23(1), 1/1997}\label{igor}

\item{P.Lancaster: `Theory of matrices', Acad. Press, New York,
    Ch. 8.4, Theo. 1, 1969}\label{lan} 

\item{D. Lenz, P. M\"uller, I. Veseli\'c: `Uniform existence of the integrated
    density of states for models on $\Z^d$', (preprint),
    arXive:math-ph/0607063v1, 2006}\label{lemuve} 

\item{C.-K. Li: `A short proof of interlacing inequalities on normalized
    Laplacians', Lin. Alg. Appl. 414, 425-427, 2006}\label{li} 

\item{R.Lyons, Y.Peres: `Probability on trees and networks', book to be
    published, C.U.P.}\label{lyons}

\item{P.Mathieu, E. Remy: `Isoperimetry and heat kernel decay on percolation
    clusters', Annals of Probability, 32 1A,  100-128, 2004}\label{matrem} 

\item{P. M\"uller, P. Stollmann: `Spectral asymptotics of the Laplacian on
    supercritical bond-percolation graphs', Journ. Func. Anal. 252, Issue 1, 1
    11 , 233-246, 2007 }\label{muesto} 

\item{D. Revelle: `Heat kernel asymptotics on the lamplighter group'.
    Electron. Comm. Probab.  8, 142--154, 2003}\label{revelle} 

\item{L.Saloff-Coste: `Lectures on finite Markov Chains', Saint-Flour,
    LNM 1665, 1997}\label{laurent}

\item{R. Schonman: `Stability of infinite clusters in supercritical
    percolation', Prob. Theory Related Fields, 113, 287-300,
    1999}\label{schon}

\item{B. Simon: `Lifshitz Tails for the Anderson Model', J. Stat. Phys., 38,
   Nrs. 1/2, 1985}\label{barry}

\item{Soardi, W.Woess: `Amenability, unimodularity, and the spectral radius of
    random walks on infinite graphs', Math. Zeitsch. 205, 471-486,
    1990}\label{sw}

\item{D.B. West: `Introduction to Graph Theory', Prentice-Hall, 1996}\label{doug} 

\item{W.Woess: `Random Walks on Infinite Graphs and Groups',
    Cambridge Tracts in Mathematics 138, Cambridge University Press, 2000
  }\label{woess}  

\end{enumerate}
\end{document}